\documentclass[a4paper, 11pt, final]{article}
\usepackage[latin1]{inputenc}
\usepackage[T1]{fontenc}     
\usepackage{fullpage}
\usepackage{amsmath}
\usepackage{dsfont}\let\mathbb\mathds
\usepackage{setspace}
\usepackage{amsfonts}
\usepackage{amssymb}
\usepackage[all,knot,poly]{xy}
\usepackage{color}
\usepackage{transparent}
\usepackage{graphicx}
\usepackage[top=3cm, bottom=3cm, left=3cm, right=3cm]{geometry}
\everymath{\displaystyle} 
\title{Asymptotic formulas for curve operators in TQFT}
\author{Renaud Detcherry \footnote{Centre de Math\'ematiques Laurent Schwartz
\'Ecole Polytechnique
Route de Saclay, 91128 Palaiseau CEDEX
France}}
\newtheorem{theo}{Theorem}[section]
\newtheorem{deftheo}[theo]{Definition and theorem}
\newtheorem{prop}[theo]{Proposition}
\newtheorem{lem}[theo]{Lemma}
\begin{document}
\maketitle
\vspace*{3 cm }
\begin{center}
\textbf{Abstract}
\newline
\newline 
\end{center}
\begin{small}
  Topological quantum field theories with gauge group $\textrm{SU}_2$ associate to each surface with marked points $\Sigma$ and each integer $r>0$, vector spaces $V_r (\Sigma)$ and to each simple closed curve $\gamma$ in $\Sigma$ Hermitian operators $T_r^{\gamma}$on the spaces $V_r (\Sigma)$. We show that the matrix elements of the operators $T_r^{\gamma}$ have an asymptotic developpement in orders of $\frac{1}{r}$, and give a formula to compute the first two terms in terms of trace functions, generalizing results of \cite{MP}.
 \newline This property, proved in \cite{MP} for the punctured torus and the 4-holed sphere, was used by them to show that curve operators on these surfaces are Toeplitz operators and to obtain the semiclassical behavior of some quantum invariants.
 \newline Our result would be a first step in order to generalize this approach to arbitrary surfaces.
 \end{small}
\newline
\newline
\newline
\textbf{Introduction}
\newline
After the discovery of the Jones polynomial, Witten defined in 1989 in \cite{W}, by a method using Feynman path integrals, a family of new invariants of 3-manifolds, together with a structure of topological quantum field theory. Reshitikhin and Turaev formalised the ideas of Witten to construct their well-known family of 3-manifolds invariants \cite{RT}: the invariants $Z_{2 r}(M), \ r \in \mathbb{N}^{*}$, together with a TQFT-structure for these invariants. An other method, more combinatorial, to define 3-manifold invariants and a TQFT from the ideas of Witten was later developped by \cite{BHMV}, this is the construction we will use in this paper, and that we will sketch in Section \ref{sec:matcoef}.
\newline The TQFT defined in \cite{BHMV} associates to each marked surface $\Sigma$ vector spaces $V_r(\Sigma)$ and to each cobordism containing a tangle $(M, T,\Sigma_0,\Sigma_1)$ a morphism $V_r(\Sigma_0) \rightarrow V_r(\Sigma_1)$. In particular, the construction associates to $\gamma$ a simple closed curve on $\Sigma$ a curve operator $T_r^{\gamma} \in \textrm{End}(V_r(\Sigma))$ (actually morphisms are associated to extended cobordisms of extended surfaces, but we can drop the extended structure when we consider only curve operators).
\newline The vector spaces $V_r(\Sigma)$ come with a natural Hermitian form. Moreover, for each banded graph $\Gamma$, such that the boundary of a regular neighboorhood of $\Gamma$ is $\Sigma$, \cite{BHMV} described a natural Hermitian basis $(\varphi_c)_{c \in U_r }$, labelled by a certain set $U_r$ of admissible coloring $c:E\rightarrow \mathbb{Z}$ of the set $E$ of edges of $\Gamma$. 
\newline Then the matrix coefficients of the curve operators admit a presentation as below:
\begin{displaymath}T_r^{\gamma}\varphi_c =\overline{c}(\gamma)\sum_{k :E \rightarrow \mathbb{Z}} F_k^{\gamma}(\frac{c}{r},\frac{1}{r})\varphi_{c+k}
\end{displaymath}
Here, $\overline{c}$ is an element of $H^{1}(\Gamma,\partial \Gamma,\mathbb{Z}/2)$ associated to an admissible color $c$, and $\overline{c}(\gamma)$ is a sign factor, always trivial if $\Gamma$ is planar.  Furthermore, we have an open set $U \subset \mathbb{R}^{E}$, such that for $c \in U_r$, we have $\frac{c}{r} \in \overline{U}$, and $F_k^{\gamma}$ are smooth functions on a neighborhood of $U \times \lbrace 0 \rbrace$ in $U \times [0,\infty)$.
\newline This result extends a result of \cite{MP}, which studied the case of curve operators on the punctured torus and the four-holed sphere. In this paper, an asymptotic formula was given to characterize the first two asymptotic terms of the $F_k^{\gamma}$ in orders of $\frac{1}{r}$ as Fourier coefficients of trace functions $f_{\gamma}$ on the modular space $\mathcal{M}(\Sigma)=\textrm{Hom}(\pi_1(\Sigma),\textrm{SU}_2)/\textrm{SU}_2$, such that $f_{\gamma} :\rho \mapsto -\textrm{tr}(\rho (\gamma))$. We generalize the asymptotic formula in \cite{MP}. Indeed, for each $q$ relative spin-structure on $(\Gamma, \partial \Gamma)$, we write $\chi(\gamma)=(-1)^{q(\gamma)}$ and we introduce:
\begin{displaymath}\sigma_{\chi}^{\gamma}(\tau ,\theta ,\hbar)=\underset{e \in E}{\sum}F_k^{\gamma}(\tau ,\hbar)e^{i k \theta}\chi(\gamma)
\end{displaymath} 
the so-called $\psi$-symbol of $\gamma$. Then we have the asymptotic expansion:
\begin{displaymath}\sigma_{\chi}^{\gamma}(\tau ,\theta ,\hbar))=f_{\gamma}(R_{\chi}(\tau ,\theta))+\frac{\hbar}{2 i}\underset{e \in E}{\sum}\frac{\partial^2}{\partial \tau_e \partial \theta_e}f_{\gamma}(R_{\chi}(\tau ,\theta)) +o(\hbar)
\end{displaymath} 
where $f_{\gamma}$ is the trace function on $\mathcal{M}(\Sigma)$ introduced above, and $R_{\chi}$ are some specific action-angle parametrizations of $\mathcal{M}(\Sigma)$. More precisely, the graph $\Gamma$ induces a pants decomposition of $\Sigma$ by curves $(C_e)_{e \in E}$ and a momentum mapping $h : \rho \rightarrow \left(\frac{1}{\pi}\textrm{Acos} (-\frac{f_{C_e}}{2})\right)_{e \in E}$, for which $\tau_e$ corresponds to action coordinates, and $\theta_e$ to angle coordinates.	
\newline
\newline The proof of \cite{MP} in the case where $\Sigma$ is the punctured torus and the four-holed sphere relied on explicit computations for some simple set of curves that generated the Kauffman algebra of $\Sigma$, then extending the result to general curves. This approach failed in higher genus as no simple set of generators is known. Instead, we developped a more conceptual and systematic method, which relies on the study of algebraic properties of the $\psi$-symbol and the Kauffman algebra of $\Sigma$.
\newline
\newline The authors of \cite{MP} used the asymptotic estimation to construct a framework for curve operators on the punctured torus and the four-holed sphere as Toeplitz operators on the sphere. This allowed to implement the WKB-approximation for eigenvectors, and deduce asymptotic expansions of quantum invariants (such as a new proof of the asymptotic expansion of $6j$-symbols, and an expression for the punctured S-matrix). Therefore, we hope to use our asymptotic expansions for general marked surface as a first step to give a framework of curve operators as Toeplitz operators on toric varieties, or at least apply the tools of microlocal analysis. Such a Toeplitz framework for curve operators would be a useful tool to study combinatorial TQFT. Indeed, in a different approach, Andersen introduced some geometrical curve operators in \cite{A06} that are Toeplitz operators to prove the asymptotic fidelity of the quantum representations of the mapping class group. We think that viewing the standard curve operators as Toeplitz operators could provide other interesting applications and hopefully a original approach to the Witten conjecture for the expansion of Reshetikhin-Turaev invariants.
\newline
\newline \textbf{Acknowledgements} We would like to thank Julien Marché for many valuable discussions.
\section{Matrix coefficients of curve operators}
\label{sec:matcoef}

\subsection{A quick overview of TQFT and curve operators}
\label{sec:TQFT}
\begin{figure}[!h]
  \begin{center}
    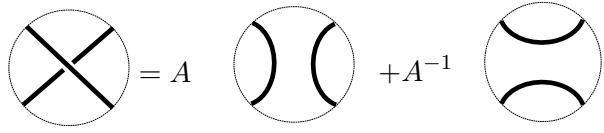
  \end{center}
  \caption{The first Kauffman relation. The other relation states that any trivial component is identified with $-A^{2}-A^{-2}$}
  \label{fig:Kauff}
  \end{figure}

To each surface $\Sigma$ with marked points $p_i$ colored by elements $\hat{c}_i$ of $\mathcal{C}_{r}=\lbrace 1, \cdots , r-1 \rbrace$, neglecting the so-called framing anomaly, the construction of \cite{BHMV} associates a vector space $V_{r}(\Sigma, \hat{c})$ and to any cobordism $(M,\Sigma_{0},\Sigma_{1})$ containing a link $L$ a morphism 
\begin{displaymath}V_{r}(M,L) \ : \ V_{r}(\Sigma_{0}) \rightarrow V_{r}(\Sigma_{1})
\end{displaymath}
 such that for every closed orientable 3-manifold M, we have $V_{r}(M) =Z_{2 r}(M)$.
\newline
\newline To give a more explicit picture of the TQFT and the vector spaces $V_{r}(\Sigma)$, we first need to introduce the notion of Kauffman bracket skein modules of 3-manifolds and Kauffman algebras of marked surfaces. 
\newline For M a 3-manifold (which can have a boundary), we define $K(M ,A)$ as the quotient of free $\mathbb{C}[A^{\pm 1}]$-module generated by links modulo isotopy and the Kauffman relations (see Figure \ref{fig:Kauff}).
\newline For $t \in \mathbb{C}^*$, we can define a Kauffman module evaluated in $t$: we write $K(M , t)=K(M , A) \underset{A=t}{\otimes}\mathbb{C}$.
\newline Now if $\Sigma$ is a surface with marked points $p_1 , \ldots , p_n$, we denote by $K(\Sigma ,A)$ the Kauffman module $K((\Sigma \setminus \lbrace p_1 , \ldots ,p_n \rbrace) \times [0,1],A)$ .
\newline We call a \textit{multicurve} on $\Sigma$ a disjoint union of simple curves on $\Sigma$, which is disjoint with the marked points of $\Sigma$. It is easy to see that $K(\Sigma ,A)$ is spanned by multicurves on $\Sigma$, and actually multicurves give a basis of this vector space. 
\newline The module $K(\Sigma ,A)$ has an algebra structure: the product $\gamma \cdot \delta$ of two elements of $K(\Sigma ,A)$ is obtained by isotopying $\gamma$ and $\delta$ so they are included in $\Sigma \times (\frac{1}{2};1]$ and $\Sigma \times [0;\frac{1}{2})$ respectively, then gluing the two parts into $\Sigma \times [0,1]$. 
\newline For $t \in \mathbb{C}^* $, we define $K(\Sigma , t)=K(\Sigma ,A) \underset{A=t} \otimes \mathbb{C}$, which is also an algebra, and admits the set of multicurves as a basis. Using this basis, we can identify $K(\Sigma , t)$ with $K(\Sigma , -1)$ and we embed $K(\Sigma , -e^{\frac{i \pi \hbar}{2}})=K(\Sigma , A) \underset{A=-e^{\frac{i \pi \hbar}{2}}}{\otimes} \mathbb{C}[[\hbar]]$ into $K(\Sigma ,-1)[[\hbar]]$. 
\newline
 \newline
The vector spaces $V_{r}(\Sigma,\hat{c})$ have a definition as quotients of Kauffman modules at roots of unity, as below:
\newline
\newline
\begin{deftheo} \label{TQFTdef} 
Let $H$ be a handlebody with $\partial H =\Sigma$ surface with marked points $p_1,\ldots ,p_n$. Given a coloration $\hat{c}$ of marked points, we chose $c_i -1$ points in a small neighboorhood of $p_i$ for each $i$, and write $P$ for the set of all such  points. We define a relative Kauffman module $K(H,\hat{c},\zeta_r)$ as the $\mathbb{C}[A^{\pm 1}]$-module generated by banded tangles in $H$ whose intersection with $\Sigma$ is the set $P$.
 For $r \in \mathbb{N}^{*}$, we write $\zeta_{r}=-e^{\frac{i \pi}{2 r}}$. For any embedding $j$ of $H$ in $\mathbf{S}^{3}$, we define the following sub-module of $K(H,\hat{c},\zeta_{r})$:
\begin{displaymath} N^{j}_{r}=\lbrace x \in K(H,\hat{c},\zeta_{r}) \ / \ \forall y \in K(\mathbf{S}^{3} \setminus \textrm{Im}(j),\hat{c},\zeta_{r}), \ \ \langle x |\underset{i=1}{\overset{r}{\otimes}} f_{c_i-1} |y \rangle =0 \rbrace 
\end{displaymath}
where we write $f_k$ for the $k$-th Jones-Wenzl idempotent, and $\langle x|\underset{i=1}{\overset{r}{\otimes}} f_{c_i-1} |y \rangle$ stands for the element of $K(\mathbf{S}^{3},\zeta_{r})$ obtained from x and y by pasting H with $\mathbf{S}^{3}\setminus \textrm{Im}(j)$, inserting Jones-Wenzl idempotent at each marked points. Then $N^{j}_{r}$ is in fact independent of $j$, is of finite codimension, and we may define: 
\begin{displaymath}V_{r}(\Sigma,\hat{c})=K(H,\hat{c},\zeta_{r})/N^{j}_{r}
\end{displaymath}
\end{deftheo}
With this setting, there is a simple description of the curve operator $T_r^{\gamma}$ associated to a multicurve $\gamma$ on $\Sigma$ disjoint from the marked points $p_1 ,\ldots ,p_n$, or more generally to an element of $K(\Sigma, \zeta_{r})$.
\newline Indeed, we can take a element $z$ of $K(H,\hat{c},\zeta_r )$ and stack a multicurve $\gamma$ over it to obtain another element $\gamma \cdot z$ of $K(H,\hat{c},\zeta_r )$. The induced map factors through $N^j_r$, as for $n \in N_r^j$ and for any $z \in K(\mathbf{S}^3 \setminus \textrm{Im}(j),\hat{c},\zeta_r)$, we have $\langle  \gamma \cdot n |\underset{i=1}{\overset{r}{\otimes}} f_{c_i-1} | z \rangle=\langle n |\underset{i=1}{\overset{r}{\otimes}} f_{c_i-1} | \gamma \cdot z \rangle$. Thus we have defined an endomorphism $T_{r}^{\gamma}$ of $V_r (\Sigma,\hat{c})$ associated to $\gamma \in K(\Sigma , \zeta_r)$.
\newline 
\newline Furthermore, the map $\begin{array}{lrcl}
T_{r}^{\cdot} : & K(\Sigma ,\zeta_r) & \longrightarrow & \textrm{End}(V_{r}(\Sigma,\hat{c})) \\
    & \gamma & \longmapsto & T_r^{\gamma} \end{array}$
 is a representation.
\newline
\newline In \cite{BHMV} a Hermitian structure on $V_{r}(\Sigma,\hat{c})$, coming from the bracket $\langle \cdot , \cdot \rangle$ that we introduced above, was provided and a Hermitian basis of $V_{r}(\Sigma,\hat{c})$ was exhibited:
\newline Let $\mathcal{C}_{r}=\lbrace 1, \ldots , r-1 \rbrace$  be the set of colors (we shifted all colors by 1 comparing to the conventions of \cite{BHMV}). We can construct a basis of the space $V_{r}(\Sigma , \hat{c})$ by the following procedure:
\newline We start by considering a pants decomposition of $\Sigma$ by a family of curves $\mathcal{C}=\lbrace C_e \rbrace_{e \in E}$ containing the components of $\partial \Sigma$. We also choose a banded graph $\Gamma$ drawn on the surface $\Sigma$
with a trivalent vertex $v_P$ lying in each pants P of the decomposition, for every $e \in E$ an edge (that we will also call $e$) joining two trivalent vertices and intersecting once the curve $C_e$ and disjoint from the other curves $C_f$, and finally $n$ univalent vertices labeled by $p_1 ,\ldots , p_n$ corresponding to the marked points of $\Sigma$. A such graph is said \textit{compatible} to the pants decomposition $\mathcal{C}$.
\newline
\newline We call admissible coloring of $\Gamma$ a map
$c : E \rightarrow C_r$ such that the following conditions hold:
\newline
- for each edge $e$ connected to a univalent vertex $p_i$ one has
 $c_e = \hat{c_i}$.
\newline
-for any triple of edges $e, f, g$ adjacent to the same vertex one has
\newline
\begin{center}
 - $c_e + c_f + c_g$ is odd
\newline - $c_e + c_f < c_g$
\newline - $c_e + c_f + c_g < 2r$.
\newline
\newline
\end{center}

We will denote by $U_r$ the set of such admissible colorations
The construction of \cite{BHMV} provides for each admissible coloring $c$ a
vector $\varphi_{c} \in V_r(\Sigma , \hat{c})$ obtained by cabling the graph $\Gamma$ by a specific combination of multicurves (we will detail this construction in Section \ref{sec:comput}). Moreover, the family $(\varphi_{c})$ when $c$ runs over all
admissible colorings is a Hermitian basis of $V_r(\Sigma , \hat{c})$.
\newline
\newline For a multicurve $\gamma$, the operators $T_r^{\gamma}$ are Hermitian operators for the Hermitian structure on $V_r(\Sigma,\hat{c})$ given by \cite{BHMV}. The spectrum and the eigenvectors of $T_{r}^{\gamma}$ are known:
\newline First, as all component of $\gamma$ are disjoint, there exists a pants decomposition of $\Sigma$ by a family of curve $\mathcal{C}=\lbrace C_{e} \rbrace_{e \in E}$ such that $\Gamma$ can be isotoped to the union of $n_e$ parallel copies of $C_e$, for some integers $n_e \in \mathbb{N}$. Then the Hermitian basis $(\varphi_c)$ coming from the pants decomposition $\mathcal{C}$ is an eigenbasis of $T_{r}^{\gamma}$, and we have 
\begin{displaymath}T_{r}^{\gamma}\varphi_c = 
\left( \underset{e \in E}{\prod}\big( -2 \cos(\frac{\pi c_{e}}{r})\big)^{n_e} \right)\varphi_c
\end{displaymath}
We should take note that the spectral radius $||T_r^{\gamma}||$ is thus always less than $2^{n(\gamma)}$, where we write $n(\gamma)$ for the number of components of the multicurve $\gamma$.
\newline
\newline In this paper, we will make a fundamental use of the following theorem that describes the Kauffman algebra $K(\Sigma ,-1)$: 
\newline 
\newline 
\begin{theo} \label{regfunc}
We note $\mathcal{M}'(\Sigma)=\textrm{Hom}(\pi_1(\Sigma) ,\textrm{SL}_{2}(\mathbb{C}))//\textrm{SL}_{2}(\mathbb{C})$ the space of characters of the fundamental group of $\Sigma \setminus \lbrace p_1 ,\cdot ,p_n \rbrace$ in $\textrm{SL}_{2}(\mathbb{C})$. This space is actually an affine algebraic variety. Let also $\textrm{Reg}(\mathcal{M}'(\Sigma))$ denote the algebra of regular functions from $\mathcal{M}'(\Sigma)$ to $\mathbb{C}$. The map 
\newline
\newline 
  $\begin{array}{lrcl} \sigma : &  K(\Sigma ,-1) & \longrightarrow & \textrm{Reg}(\mathcal{M}'(\Sigma)) \\
& \gamma &\longrightarrow & f_{\gamma} \ \textrm{such that } \ f_{\gamma}(\rho)=-\textrm{Tr}(\rho(\gamma))
\end{array}$
\newline
\newline
is a injective morphism of algebras.
\end{theo}
 Bullock \cite{Bul97} and  Brumfiel-Hilden \cite{BH} first independently proved this gives an isomorphism from $K(\Sigma , -1)$ modulo its nilradical to $\mathcal{M}'(\Sigma)$. Their work was later completed by Sikora-Przytycki \cite{PS00} and independently by Charles and Marché \cite{CM09} to give Theorem \ref{regfunc}.
\newline 
\newline Finally, we end this preliminary section with a formula due to Goldman for products of elements of the Kauffman algebra at first order. We recall that $\mathcal{M}'(\Sigma)$ is a Poisson manifold for the Poisson structure given in \cite{G86}. This Poisson structure depends of a choice of normalization of the symplectic structure on $\mathcal{M}(\Sigma)$. We normalize the symplectic form $\omega$ as the symplectic reduction of the form $\omega (\alpha ,\beta) =\frac{1}{2 \pi}\int_{\Sigma}\textrm{Tr}(\alpha \wedge \beta)$ for $\alpha ,\beta \in \Omega^1(\Sigma ,\textrm{su}_2)$. As by the previous theorem it is possible to link the product of elements of $K(\Sigma ,-1)$ with products of trace functions on $\mathcal{M}'(\Sigma)$, the works of Goldman \cite{G86} and Turaev \cite{TU91} gave a way to think of the first order of a product of elements in $K(\Sigma ,-e^{\frac{i \pi \hbar}{2}})$ as a Poisson bracket of trace functions.
\newline
\newline
\begin{theo} \label{prodform}
Let $\gamma$ and $\delta$ be two multicurves, which can be viewed as elements of $K(\Sigma ,-e^{\frac{i \pi \hbar}{2}})=K(\Sigma ,A) \underset{A=-e^{\frac{i \pi \hbar}{2}}}{\otimes} \mathbb{C}[[\hbar]]$. This space is isomorphic to a subspace of $\textrm{Reg}(\mathcal{M}'(\Sigma))[[\hbar]]$ via the map $\sigma$ of Theorem \ref{regfunc}. Then we have:
\begin{displaymath} \gamma \cdot \delta =f_{\gamma} f_{\delta} + \frac{\hbar}{i} \lbrace f_{\gamma} , f_{\delta} \rbrace +o(\hbar)
\end{displaymath}
\end{theo}

\subsection{An asymptotic expression for matrix coefficients of curve operators}
\label{sec:asympt}
\begin{figure}
   \centering
  \def\svgwidth{200pt}
  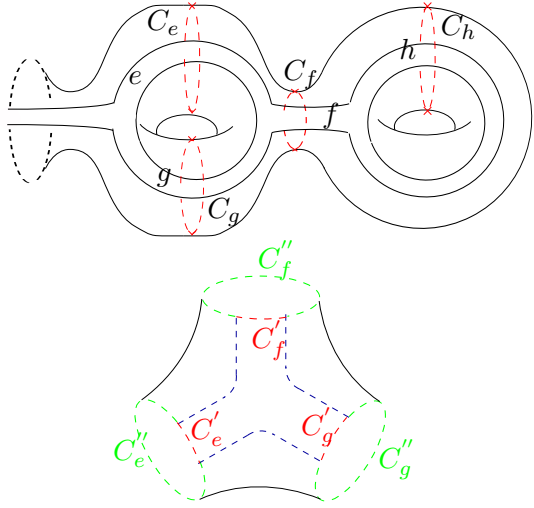
  \caption{A banded graph compatible to a pants decomposition of $\Sigma$ by the curve $C_e$ and the associated cell decomposition of a pants into hexagons}
  \label{fig:celldecomp}
\end{figure}

In this section, we fix a surface $\Sigma$ with marked points $p_1 ,\ldots , p_n$. We fix $t_1 , \ldots ,t_n \in \mathbb{Q} \cap [0,1]$. and we define colorations of the marked points as $(\hat{c}_r )_i=r t_i$. Thus we will need to suppose that $r$ is a multiple of $D$, the common denominator of the $t_i$, so that $(\hat{c}_r )_i$ are integers.
 For any multicurve $\gamma$, and $r$ a mutiple of $D$, we recall that we call $T_{r}^{\gamma}$ the endomorphism of $V_{r}(\Sigma ,\hat{c}_r)$ associated to $\gamma$ by the TQFT. We fix a pants decomposition $\lbrace C_{e} \rbrace_{e \in E}$ of $\Sigma$ and a compatible trivalent banded graph $\Gamma$ drawn on $\Sigma$, and thus a Hermitian basis $( \varphi_c)_{c \in U_{r}}$ of $V_r (\Sigma,\hat{c}_r)$ as explained above.
\newline
\newline Note that the data of $\Gamma$ give us a cell decomposition of $\Sigma$ into a bunch of hexagons, by the boundary components of $\Gamma$, and the curves $C_e$. For each $e \in E$, we name $C_e^{'}$ (respectively $C_e^{''}$) the segment $\Gamma \cap C_e$ (resp. $C_e \setminus \textrm{Int}(C_e \cap \Gamma)$), see Figure \ref{fig:celldecomp}. 
\newline We will associate to each admissible color $c$ an element $\overline{c}$ of $H^{1}(\Sigma ,\mathbb{Z}/2)$ by writing, for $\gamma \in H_1 (\Sigma ,\mathbb{Z}/2)$:
\begin{displaymath} \overline{c}(\gamma) = \underset{e \in E}{\prod} (-1)^{(c_e-1)( C_e^{' *}(\gamma)+ C_e^{'' *}(\gamma))}
\end{displaymath}
In this formula, $C_e^{' *}$ (resp. $C_e^{'' *}$) is the cellular cochain dual to $C_e^{'}$ (resp. $C_e^{''}$). We can check that the such defined $\overline{c}$ is indeed a cocycle as its value on the boundary of each hexagon is of the form $(-1)^{c_e +c_f+c_g-1}$ for $e,f,g$ three adjacent edges, which equals to $1$ as $c$ is an admissible color.
\newline
We introduce a set U of "real admissible colorations", that we view as a set of limits of admissible colorations in $U_r$. It consists of maps $E \rightarrow (0,1)$ such that for any triple $e$,$f$, and $g$ of edges of the graph associated to the pants decomposition that are adjacent to the same trivalent vertex, we have:
\newline
\begin{center}
\begin{itemize} \item[-] $\tau_{e}+\tau_{f} < \tau_{g}$
 \item[-] $\tau_e+\tau_f+\tau_g < 2$ 
 \item[-] for $e \in E$ adjacent to a marked point $p_i$ we have $\tau_e=t_i$.
\newline
\end{itemize}
 \end{center}
We can now give an expression of $T_{r}^{\gamma} \varphi_{c}$ when $\frac{c}{r}$ tends to an element $\tau \in U$ when $r \rightarrow + \infty$.
\newline
\newline
\begin{theo} \label{maintheo1}
Let $\gamma$ be a multicurve in $\Sigma \setminus \lbrace p_1 , \ldots , p_n \rbrace$
\newline 1) There is an open set $V_{\gamma} \subset U \times [0,1]$ containing $U \times \lbrace 0 \rbrace$ and functions $(F_{k}^{\gamma})_{k : E \rightarrow \mathbb{Z}}$ that are smooth on $V_{\gamma}$ such that we have for any $c \in U_{r}$,
\begin{displaymath}T_{r}^{\gamma} \varphi_{c} = \overline{c}(\gamma) \underset{k : E \rightarrow \mathbb{Z}}{\sum} F_{k}^{\gamma}(\frac{c}{r},\frac{1}{r}) \varphi_{c+k}
\end{displaymath}
\newline
2)Let $I_e= \sharp (\gamma \cap C_e)$. If there exists $e \in E$ such that $k_{e} >I_{e}$ or such that $k_e \neq I_e (\textrm{mod} \ 2)$, then $F_{k}^{\gamma}=0$.
\end{theo}
We remind that here, $\overline{c}$ is an element of $H^1 (\Sigma , \mathbb{Z}/2 )$, so that $\overline{c}(\gamma)$ is just a sign. So up to this oscillating sign which varies with the level $r$, the matrix coefficients of curves operators are converging when $r \rightarrow +\infty$. One might remark that this sign factor did not appear in \cite{MP}, but it can be shown that it is trivial when the banded trivalent graph $\Gamma$ is planar (which was the case for the punctured torus and the four-holed sphere).
\newline The two points of the theorem are proved by doing local computations, working in each pants of the pants decomposition of $\Sigma$. The proof is rather technical but not difficult and relies on fusion rules in TQFT. It will be detailed in \ref{sec:comput}.
\newline The open sets $V_{\gamma}$ are actually explicit: we will show that we can take 
\newline $V_{\gamma}=\lbrace (\tau ,\hbar ) \ / \ (\tau_e + \varepsilon_e \hbar I_e)_{e \in E} \in U \ , \ \forall \varepsilon \in \lbrace \pm 1 \rbrace ^{E} \rbrace$ 
\newline
\newline The coefficients $F_k$ can be computed recursively for any surface $\Sigma$ (together with pants decomposition $\mathcal{C}$ and trivalent banded graph $\Gamma$), but are uneasy to make explicit. However, we will provide a formula to get the first two term $F_k^{0}$ and $F_k^{1}$ of their asymptotic expansion in orders of $\hbar=\frac{1}{r}$: 
\begin{displaymath}F_k(\tau ,\hbar)=F_k^{0}(\tau) +\hbar F_k^{1}(\tau) + o(\hbar)
\end{displaymath}
To state our asymptotic formula, we have yet another definition to give.
\newline 
\newline For $\Gamma$ a trivalent banded graph on $\Sigma$ compatible to the pants decomposition $\mathcal{C}$, we define the \textit{intersection algebra} $A_{\Gamma}$ as
\begin{displaymath}A_{\Gamma}=\underset{\alpha \in H_1(\Gamma , \partial \Gamma ,\mathbb{Z}/2)}{\bigoplus}\mathbb{C} \cdot [\alpha]
\end{displaymath}
with the product structure given by $[\gamma] \cdot [\delta] = \langle \gamma , \delta \rangle [\gamma +\delta]$. Here $\langle \cdot , \cdot \rangle$ shall be understood as the intersection form in $H_1 (\Gamma , \partial \Gamma , \mathbb{Z}/2)$,  The associativity of the product in $A_{\Gamma}$ follows easily from the bilinearity of the intersection form.
\newline Note that $A_{\Gamma}$ is a commutative (the intersection form in $H_1 (\Gamma , \partial \Gamma , \mathbb{Z}/2)$ is symmetric) $\mathbb{C}$-algebra of dimension $2^{g}$, where $g$ is the genus of $\Sigma$. Define $\hat{A}_{\Gamma}$ the set of morphisms $A_{\Gamma} \rightarrow \mathbb{C}$, then $A_{\Gamma}$ is isomorphic to $\mathbb{C}^{\hat{A}_{\Gamma}}$, by the map $a \in A_{\Gamma} \mapsto (\chi (a))_{\chi \in \hat{A}_{\Gamma}}$. 
\newline
Remember that $\Gamma$ and $\mathcal{C}$ give a cell decomposition of $\Sigma$ into hexagons, with two hexagons associated to each trivalent vertex of $\Gamma$. In fact, we have a continuous map $p : \Sigma \rightarrow \Gamma$ obtained by identifying the hexagons associated to the same vertex, and we write $p_*$ for the map $H_1 (\Sigma ,\mathbb{Z}/2) \rightarrow H_1 (\Gamma ,\partial \Gamma , \mathbb{Z}/2)$ induced by $p$. 
\newline
We can now turn to the statement of our asymptotic formula for the matrix coefficients $F_k$ of curves operators:
\newline
\newline  
\begin{theo} \label{maintheo2}
Let $\gamma$ be a multicurve on $\Sigma$. For all $\theta \in (\mathbb{R}/2 \pi \mathbb{Z})^{E}$ and $(\tau , \hbar) \in V_{\gamma}$,
\newline we define the $\psi$-symbol of $T_r^{\gamma}$ as the following element of $A_{\Gamma}$:
\begin{displaymath}\sigma^{\gamma}(\tau , \theta , \hbar )=\underset{k:E \rightarrow \mathbb{Z}}{\sum} F_{k}(\tau , \hbar) e^{i k \cdot \theta} [p_* (\gamma)]
\end{displaymath} 
Then, $\sigma^{\gamma}(\tau , \theta , \hbar)$ has the following asymptotic expansion:
\begin{displaymath}\sigma^{\gamma}(\tau ,\theta ,\hbar ) = \sigma^{\gamma}(\tau ,\theta ,0) + \frac{\hbar}{2 i} \underset{e \in E}{\sum}\frac{\partial^{2}}{\partial \tau_e \partial \theta_e} \sigma^{\gamma}(\tau , \theta ,0) +o(\hbar)
\end{displaymath}
and writing $\sigma_{\chi}^{\gamma} (\tau ,\theta) = \chi (\sigma^{\gamma} (\tau ,\theta ,0))$, we have $\sigma_{\chi}^{\gamma} (\tau ,\theta)= f_{\gamma}(R_{\chi} (\tau ,\theta))= -\textrm{Tr}(R_{\chi} (\tau ,\theta)(\gamma))$, where the $R_{\chi}$ are action-angle parametrizations on $\mathcal{M}(\Sigma)=\textrm{Hom}(\pi_1 (\Sigma \setminus \lbrace p_1 ,\ldots ,p_n \rbrace ) , \textrm{SU}_2)/\textrm{SU}_2$, differing only by the origin of angles.
\end{theo}
More precisely, there is a moment mapping on $\mathcal{M}(\Sigma)$ given by
\begin{displaymath}
 \begin{array}{lrcl} h : &\mathcal{M}(\Sigma) & \mapsto & U \\
                    & \rho & \mapsto &  (h_{C_e}(\rho)=\tau_e (\rho)=\frac{1}{\pi} \textrm{Acos}(\frac{1}{2} \textrm{Tr}(\rho (C_e))))_{e \in E}
                    \end{array}
\end{displaymath}
The variables $h_{C_e}$ are independent Poisson commuting functions on $\mathcal{M}(\Sigma)$ and their Hamiltonian flow gives an action of a torus on $\mathcal{M}(\Sigma)$, and thus induces \textit{angle coordinates} $\theta_e$ (unique up to change of origin) on each level set of the $h_{C_e}$. The representation $R_{\chi}(\tau ,\theta)$ has $(\tau_e ,\theta_e)$ as coordinates on $\mathcal{M}_{\textrm{irr}}(\Sigma)$.
\newline
\newline
Let us add a few remarks on the definition of the $\psi$-symbol:
\newline 1) To begin with, the sum over $k : E \rightarrow \mathbb{Z}$ is actually a finite sum, as only a finite number of coefficients $F_{k}^{\gamma}$ does not vanish by the second point of Theorem \ref{maintheo1}. Furthermore, we wrote $k \cdot \theta$ for $\underset{e \in E}{\sum} k_{e} \theta_{e}$. Also, we will often omit the $p_*$ and just write $[\gamma]$ for the element $[p_*(\gamma)]$, when $\gamma$ is a multicurve.
\newline 2) We can recover the matrix coefficents $F_{k}^{\gamma}$ from the $\psi$-symbol by taking Fourier coefficients of $\sigma^{\gamma}(\tau , \cdot , \hbar )$.
\newline 3) We will often refer to the zero order in $\hbar$ of the $\psi$-symbol, that is $\sigma^{\gamma}(\tau , \theta ,0)$, as the principal symbol of $T_{r}^{\gamma}$.
\newline 4) For a fixed $(\tau , \theta , \hbar )$, this definition only introduce $\gamma \mapsto \sigma^{\gamma}(\tau , \theta , \hbar )$ as a map from multicurves to $A_{\Gamma}$. We extend it by multilinearity to obtain $\sigma(\tau , \theta , \hbar ) \ : \ K(\Sigma, -e^{\frac{i \pi \hbar}{2}}) \rightarrow A_{\Gamma}[[\hbar]]$, as $K(\Sigma, -e^{\frac{i \pi \hbar}{2}})$ is spanned by multicurves.
\newline
\newline
The proof of Theorem \ref{maintheo2}, giving an asymptotic formula for the $\psi$-symbol, will be the goal of Sections \ref{sec:prsym} and \ref{sec:frstor}. It will rely heavily on the following property of the $\psi$-symbol, that explains its compatibility with the product in $K(\Sigma ,-e^{\frac{i \pi \hbar}{2 }})$:
\newline
\newline
\begin{prop}\label{mult}
Let $\gamma$ and $\delta$ be two multicurves on $\Sigma$. Then we have the following asymptotic expression:
\begin{displaymath} \sigma^{\gamma \cdot \delta}(\tau , \theta , \hbar ) =\left(\sigma^{\gamma}(\tau , \theta , \hbar) \sigma^{\delta}(\tau , \theta , \hbar ) + \frac{\hbar}{i} \underset{e}{\sum} \partial_{\tau_{e}}\sigma^{\gamma}(\tau , \theta , \hbar ) \partial_{\theta_{e}} \sigma^{\delta}(\tau , \theta , \hbar ) \right) + o(\hbar) 
\end{displaymath} 
\end{prop}
According to this proposition, the principal symbol $\sigma^{\cdot}(\tau , \theta , 0) \ : \ K(\Sigma ,-1) \rightarrow A_{\Gamma}$ is a morphism of algebras. Recall that there is an isomorphism betweem the algebras $A_{\Gamma}$ and $\mathbb{C}^{\hat{A}_{\Gamma}}$, the components of which are morphisms we named $\chi : A_{\Gamma} \rightarrow \mathbb{C}$. Then the maps $\sigma_{\chi}=\chi \circ \sigma(\tau , \theta ,0) : K(\Sigma ,-1) \rightarrow \mathbb{C}$ constitutes a collection of algebra morphisms.
\newline We will then use the theorem of Bullock to show that these morphisms have the form $f \in \textrm{Reg}(\mathcal{M}'(\Sigma)) \rightarrow f(R_{\chi})$, for some representations $R_{\chi}$ of $\pi_1 (\Sigma \setminus \lbrace p_1 , \ldots p_n \rbrace)$. 
\newline Finally, we identify precisely the representations $R_{\chi}$ and how they depend on $(\tau ,\theta)$ by checking the special values of the $\psi$-symbol on the curves $C_e$.
\newline As for the computation of the first order term, we proceed in a similar fashion: first we will show, using only Proposition \ref{mult} that this term is related to derivations of algebras $K(\Sigma , -1) \rightarrow A$, then by studying the values of the $\psi$-symbol on the curves $C_e$ and an other family of curves $D_e$, we show the first order term is indeed given by the formula in Theorem \ref{maintheo2}.
\section{Computations of curve operators using fusion rules}
\label{sec:comput}
This section is devoted to the skein theory computations that will be
needed in order to prove Theorem \ref{maintheo1}. We describe the
general form of matrix coefficents of curve operators, and give examples
of explicit computations of the coefficients $F_k^{\gamma}$ and the
$\psi$-symbol $\sigma^{\gamma}$ for some curves $\gamma$.
\subsection{Fusion rules in a pants decomposition}
\label{sec:fusion}
\begin{figure}

  \centering
  \def\svgwidth{150pt}
  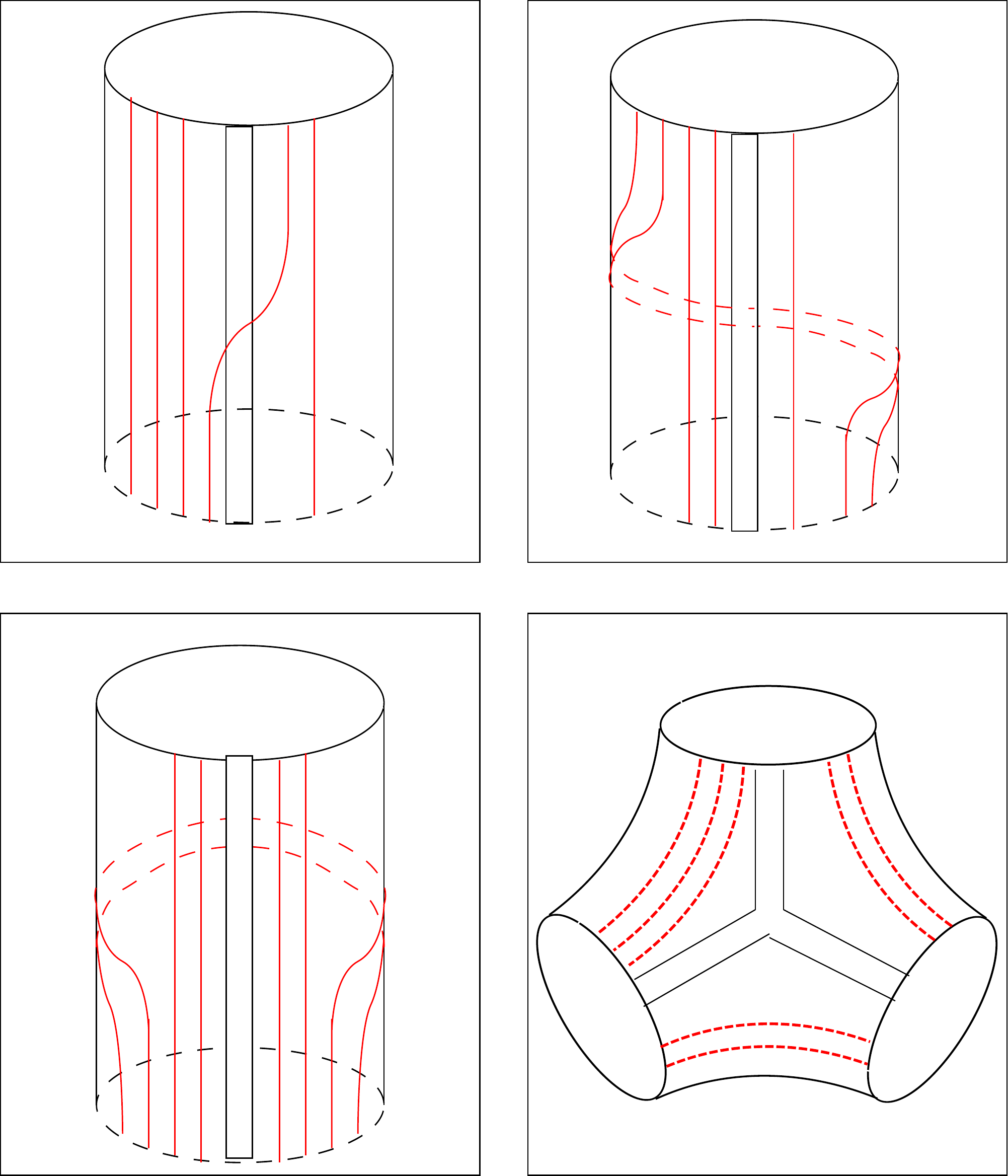
  \caption{Dehn presentation of multicurves}
  \label{fig:dehncoor}
\end{figure}

\begin{figure}[htb]
  \begin{center}
  \def\svgwidth{135pt}
  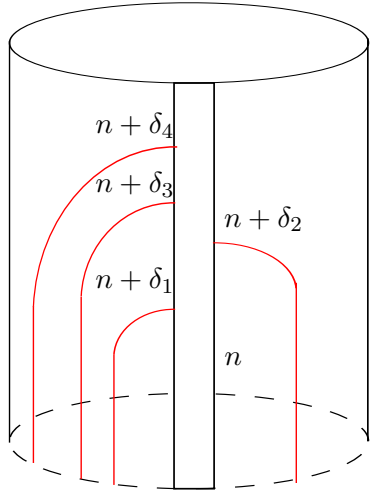
  \caption{A candlestick $C(n,\varepsilon ,\theta)$ with 4 legs }
  \label{fig:chandelier}
  \end{center}
  We wrote $\delta_i =\underset{j=1}{\overset{i}{\sum}} \varepsilon_j$ the partial sums of the shifts. In the bottom the central edge is colored by $n$, and the color is shifted by $\varepsilon_i$ when $\Gamma$ meets the i-th leg. Notice that the legs can go alternatively to the left or to the right of the central edge
  \end{figure}

In this paragraph, we will  work with a fixed surface $\Sigma$, along with a pants
decomposition by a family of curves $\mathcal{C}=\lbrace C_e
\rbrace_{e \in E}$. We can consider $n_e \geq 1$ parallel copies $(C_e^{k})_{1 \leq k \leq n_e}$ of the curves $C_e$ so that the curves
 $C_e^{k}$ cut the surface $\Sigma$ into a collection of pants
$\lbrace P_s \rbrace_{s \in S}$ and annuli $\lbrace A_e^{k} \ , e \in E \ , \ \ 1 \leq k \leq n_e -1 \rbrace$.
\newline We recall that to this pants decomposition is associated a
Hermitian basis $\varphi_c$ of $V_r (\Sigma)$, of which we will remind the construction:
\newline
\newline Let $\Gamma$ be a banded trivalent
graph compatible to the pants decomposition $\mathcal{C}$ of $\Sigma$ as in Subsection \ref{sec:TQFT}. We recall that $\Gamma$ is viewed as drawn on $\Sigma$.  Given an admissible coloring
$c \ : \ E \rightarrow C_r$, we define  $\psi_c \in  K(\Sigma; \hat{c}; \zeta_r)$ w:
\begin{itemize}
\item[-]
Replace each edge e of $\Gamma$ by $c_e -1$ parallel copies of e lying on $\Sigma$.
\item[-] 
Insert in the middle of each edge the idempotent $f_{c_e -1}$ where we write $f_k$ for the $k$-th Jones-Wenzl idempotent (see \cite{BHMV} for details).
\item[-]
In the neighborhood of each trivalent vertex, join the three bunches of
lines in $\Sigma$ in the unique possible way avoiding crossings.
\end{itemize}
This family of vectors is actually an orthogonal basis of $V_r (\Sigma , c)$ for
a natural Hermitian structure defined in \cite{BHMV}, that we do not recall here. We refer to Theorem
4.11 in \cite{BHMV} for the proof and the following formula:
\begin{equation} ||\psi_c ||^{2}=\big( \frac{2}{r}\big)^{\frac{\chi (\Gamma)}{2}} \frac{\prod_P \langle c_P^{1} , c_P^{2} , c_P^{3} \rangle}{\prod_e \langle c_e \rangle}
\end{equation}
Here the first product is over all vertex $P$ corresponding to pants of the pants decomposition, the second over the edges e of the graph $\Gamma$. We write $\langle n \rangle$ for $\sin (\frac{\pi n}{r})$, $\langle n \rangle !$ for $\prod_{i=1}^{n} \langle i \rangle$, $c_P^{1}$, $c_P^{2}$, and $c_P^{3}$ for the colors of the 3 edges adjacent to P, and we also set
\begin{displaymath}\langle a , b , c \rangle =\frac{ \langle \frac{a+b+c-1}{2} \rangle ! \langle \frac{a+b-c-1}{2} \rangle ! \langle \frac{a-b+c-1}{2} \rangle ! \langle \frac{b+c-a-1}{2} \rangle !}{\langle a-1 \rangle ! \langle b-1 \rangle ! \langle c-1 \rangle !}
\end{displaymath}
As we will work with TQFT vectors locally, inside a pants of the pants decomposition for example, we will need to give a local version of this norm. Notice that if we forget the global factor $(\frac{2}{r})^{\frac{\chi ( \Gamma )}{2}}$ in the norm, we will not change the matrix coefficients of curve operators $T_r^{\gamma}$. Then, we will decide that the square of the norm of a trivalent graph is 
\begin{displaymath}\frac{\prod_P \langle c_P^{1} , c_P^{2} , c_P^{3} \rangle}{\prod_{e \in E_2} \langle c_{e} \rangle \prod_{e \in E_1} \langle c_e \rangle^{\frac{1}{2}}}
\end{displaymath}
where the products in the denominator are over $E_2$, the set of edges adjacent to 2 internal vertices (we include the marked points here),  and $E_1$ the set of edges adjacent to 1 internal vertex and 1 external vertex (the other edges bear no contribution to the norm). With this definition, if we paste pieces of colored graph to get the graph $\Gamma$, we obtain the previous norm.
\newline
 With this setting, we give a normalized version of fusion rules in TQFT. The fusion
rules derived in \cite{MV}, give a way to compute the image of the vector
$\varphi_c$ by curves operators. We list the fusion rules that we will
need below; our version differs from the rules in \cite{MV}, as
we express them with the normalized vectors $\varphi_c$ instead of the
vectors $\psi_c$.
\begin{figure}[h!]
  \begin{center}
  \def\svgwidth{300pt}
  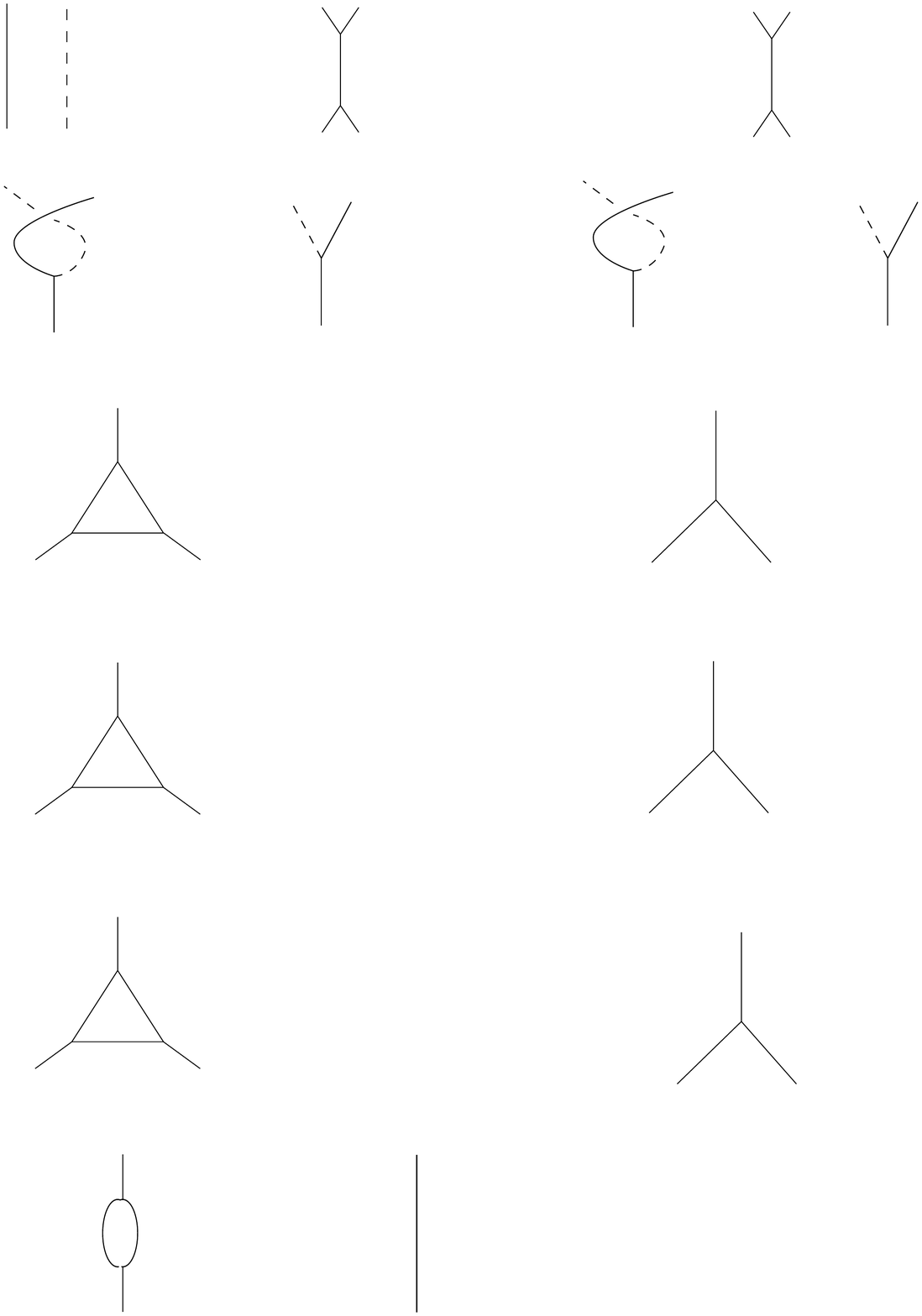
  \caption{Fusion rules}
  \label{fig:fusion}
  \end{center}
  These "normalized" fusion rules allow to simplify the union of a colored banded graph and a curve colored by 2. All edges for which we did not write a color are colored by 2. The first rule allows to merge an edge colored by 2 with an other one. The second line consists of the "half-twist formulas" of \cite{MV}. When all curves have been merged with the graph, the 3rd, 4th and 5th lines can be used to remove trigons, and the last rule to remove bigons.
  \end{figure}
\newpage
 We will lead the computations by using fusion rules only locally,
that is only inside of a pair of pants of the pants decomposition , or inside an
annulus  in the neighboorhood of one of the curves $C_e$.
\newline Indeed, for $\gamma$ a multicurve, by a classification provided
by Dehn, we can isotope $\gamma$ so that the intersection of $\gamma$ with
each pants $P_s$ of the decomposition looks like the 4th picture of Figure \ref{fig:dehncoor}, and the intersection with each of the
annulus $A_e^{k}$ looks like one of the first three pictures of Figure \ref{fig:dehncoor}.
\newline Furthermore, in this isotopy class, the intersection of $\gamma$ with each $C_e$ is the smallest in the isotopy class of $\gamma$, see \cite{FLP}.
\newline
\newline Now, we do the computations in two steps:
\newline First, we use fusion rules to reduce each type of piece to elements corresponding to the trace of the graph $\Gamma$ in a pants or annulus with a certain coloring, glued with "candlesticks". 
\newline We call a \textit{candlestick} an element of the TQFT of an annulus that is the normalized vector associated to a banded trivalent graph in an annulus, consisting of a central edge joining the boundary components (with no twist), colored by $n \in \mathcal{C}_r$ on the bottom component, a collection of legs colored by 2, joining the central edge and the bottom component, as in Figure \ref{fig:chandelier}.
 The data that define a candlestick with $k$ legs $C(n,\varepsilon, \Theta)$ is the color $n \in \mathcal{C}_r$ of the central edge at the bottom, the order $\Theta$ in which the legs join the central edge, and the shifts  of the color of the central edge $(\varepsilon_i)_{i=1 \ldots k}$ when we pass each vertex corresponding to a leg.
  \newline
  \newline
\newline \textbf{Reduction of the different pieces :}
Simple computations using fusion rules give us the following formulas when the pants or the annulus contain only one curve:
\newline
\newline
\begingroup%
  \makeatletter%
  \providecommand\color[2][]{%
    \errmessage{(Inkscape) Color is used for the text in Inkscape, but the package 'color.sty' is not loaded}%
    \renewcommand\color[2][]{}%
  }%
  \providecommand\transparent[1]{%
    \errmessage{(Inkscape) Transparency is used (non-zero) for the text in Inkscape, but the package 'transparent.sty' is not loaded}%
    \renewcommand\transparent[1]{}%
  }%
  \providecommand\rotatebox[2]{#2}%
  \ifx\svgwidth\undefined%
    \setlength{\unitlength}{318.8648969bp}%
    \ifx\svgscale\undefined%
      \relax%
    \else%
      \setlength{\unitlength}{\unitlength * \real{\svgscale}}%
    \fi%
  \else%
    \setlength{\unitlength}{\svgwidth}%
  \fi%
  \global\let\svgwidth\undefined%
  \global\let\svgscale\undefined%
  \makeatother%
  \begin{picture}(1,0.22780321)%
    \put(0,0){\includegraphics[width=\unitlength]{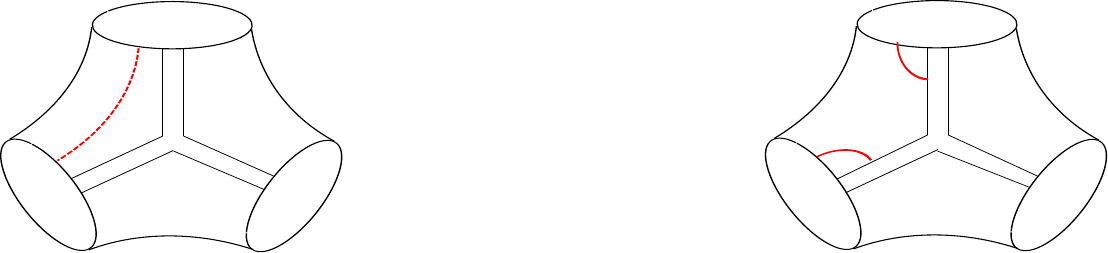}}%
    \put(0.19291302,0.14964371){\color[rgb]{0,0,0}\makebox(0,0)[lb]{\smash{$b$}}}%
    \put(0.11167429,0.03442325){\color[rgb]{0,0,0}\makebox(0,0)[lb]{\smash{$c$}}}%
    \put(0.31477111,0.03442325){\color[rgb]{0,0,0}\makebox(0,0)[lb]{\smash{$a$}}}%
    \put(0.87179555,0.19675428){\color[rgb]{0,0,0}\makebox(0,0)[lb]{\smash{$b$}}}%
    \put(0.68034565,0.03442325){\color[rgb]{0,0,0}\makebox(0,0)[lb]{\smash{$c$}}}%
    \put(0.92406175,0.02290123){\color[rgb]{0,0,0}\makebox(0,0)[lb]{\smash{$a$}}}%
    \put(0.87179555,0.10894281){\color[rgb]{0,0,0}\makebox(0,0)[lb]{\smash{$b+ \varepsilon$}}}%
    \put(0.81574339,0.03442325){\color[rgb]{0,0,0}\makebox(0,0)[lb]{\smash{$c+ \mu$}}}%
    \put(0.35747121,0.13403179){\color[rgb]{0,0,0}\makebox(0,0)[lb]{\smash{$=\underset{\varepsilon ,\mu}{\sum}F_{\varepsilon ,\mu}(a,b,c,r)$}}}%
  \end{picture}%
\endgroup%

\newline
\newline
where we set $F_{+ , +}(a,b,c,r)=\left(\frac{\langle \frac{a+b+c+1}{2} \rangle \langle \frac{b+c-a+1}{2} \rangle}{\langle b \rangle \langle c \rangle}\right)^{\frac{1}{2}}$, 
\newline $F_{+,-}(a,b,c,r)=F_{-,+}(a,c,b,r)=-\left( \frac{\langle \frac{a-b+c-1}{2} \rangle \langle \frac{a+b-c-1}{2} \rangle}{\langle b \rangle \langle c \rangle} \right)^{\frac{1}{2}}$
\newline and $F_{-,-}(a,b,c,r)=-\left(\frac{\langle \frac{a+b+c-1}{2} \rangle \langle \frac{b+c-a-1}{2} \rangle}{\langle b \rangle \langle c \rangle}\right)^{\frac{1}{2}}$
\newline 
\newline
\begin{center}
\begingroup%
  \makeatletter%
  \providecommand\color[2][]{%
    \errmessage{(Inkscape) Color is used for the text in Inkscape, but the package 'color.sty' is not loaded}%
    \renewcommand\color[2][]{}%
  }%
  \providecommand\transparent[1]{%
    \errmessage{(Inkscape) Transparency is used (non-zero) for the text in Inkscape, but the package 'transparent.sty' is not loaded}%
    \renewcommand\transparent[1]{}%
  }%
  \providecommand\rotatebox[2]{#2}%
  \ifx\svgwidth\undefined%
    \setlength{\unitlength}{298.2413514bp}%
    \ifx\svgscale\undefined%
      \relax%
    \else%
      \setlength{\unitlength}{\unitlength * \real{\svgscale}}%
    \fi%
  \else%
    \setlength{\unitlength}{\svgwidth}%
  \fi%
  \global\let\svgwidth\undefined%
  \global\let\svgscale\undefined%
  \makeatother%
  \begin{picture}(1,0.273468)%
    \put(0,0){\includegraphics[width=\unitlength]{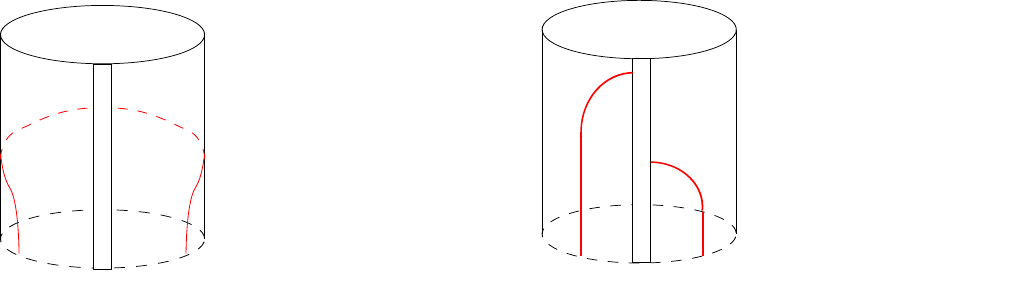}}%
    \put(0.21504589,0.0973954){\color[rgb]{0,0,0}\makebox(0,0)[lb]{\smash{$=\underset{\varepsilon}{\sum}G_{\varepsilon}(n,r)$}}}%
    \put(0.12262619,0.12165541){\color[rgb]{0,0,0}\makebox(0,0)[lb]{\smash{$n$}}}%
    \put(0.62849126,0.25116092){\color[rgb]{0,0,0}\makebox(0,0)[lb]{\smash{$n$}}}%
    \put(0.60604048,0.00191912){\color[rgb]{0,0,0}\makebox(0,0)[lb]{\smash{$n$}}}%
    \put(0.65621575,0.14571242){\color[rgb]{0,0,0}\makebox(0,0)[lb]{\smash{$n+\varepsilon$}}}%
  \end{picture}%
\endgroup%

 \end{center}
 where $G_{+}(n,r)=(-1)^{n+1}e^{\frac{-i \pi (n-1)}{r}} \left(\frac{\langle n+1 \rangle }{\langle n \rangle}\right)^{\frac{1}{2}}$
\newline and $G_{-}(n,r)=(-1)^{n+1}e^{\frac{i \pi (n+1)}{r}} \left(\frac{\langle n-1 \rangle }{\langle n \rangle}\right)^{\frac{1}{2}}$
\newline
\newline 
\begin{center}
\begingroup%
  \makeatletter%
  \providecommand\color[2][]{%
    \errmessage{(Inkscape) Color is used for the text in Inkscape, but the package 'color.sty' is not loaded}%
    \renewcommand\color[2][]{}%
  }%
  \providecommand\transparent[1]{%
    \errmessage{(Inkscape) Transparency is used (non-zero) for the text in Inkscape, but the package 'transparent.sty' is not loaded}%
    \renewcommand\transparent[1]{}%
  }%
  \providecommand\rotatebox[2]{#2}%
  \ifx\svgwidth\undefined%
    \setlength{\unitlength}{278.46674708bp}%
    \ifx\svgscale\undefined%
      \relax%
    \else%
      \setlength{\unitlength}{\unitlength * \real{\svgscale}}%
    \fi%
  \else%
    \setlength{\unitlength}{\svgwidth}%
  \fi%
  \global\let\svgwidth\undefined%
  \global\let\svgscale\undefined%
  \makeatother%
  \begin{picture}(1,0.2923349)%
    \put(0,0){\includegraphics[width=\unitlength]{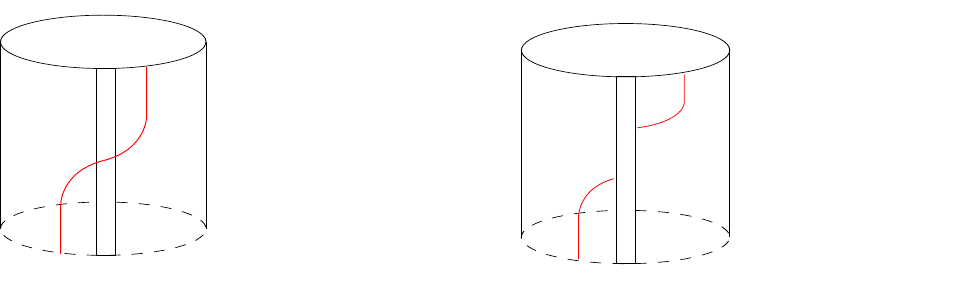}}%
    \put(0.03596082,0.1938773){\color[rgb]{0,0,0}\makebox(0,0)[lb]{\smash{$n$}}}%
    \put(0.62267178,0.27603748){\color[rgb]{0,0,0}\makebox(0,0)[lb]{\smash{$n$}}}%
    \put(0.67130405,0.0021092){\color[rgb]{0,0,0}\makebox(0,0)[lb]{\smash{$n$}}}%
    \put(0.66799318,0.13215077){\color[rgb]{0,0,0}\makebox(0,0)[lb]{\smash{$n+\varepsilon$}}}%
    \put(0.23706203,0.16087952){\color[rgb]{0,0,0}\makebox(0,0)[lb]{\smash{$=\underset{\varepsilon}{\sum}H_{\varepsilon}(n,r)$}}}%
  \end{picture}%
\endgroup%

\end{center}
 where $H_+ (n,r)=(-1)^{n+1} e^{\frac{i \pi(n-1)}{r}} \left(\frac{\langle n+1 \rangle }{\langle n \rangle}\right)^{\frac{1}{2}}$
\newline and $H_- (n,r)=(-1)^{n+1} e^{\frac{-i \pi (n+1)}{r}} \left(\frac{\langle n-1 \rangle }{\langle n \rangle}\right)^{\frac{1}{2}}$
\newline 
\newline 
\begin{center}
\begingroup%
  \makeatletter%
  \providecommand\color[2][]{%
    \errmessage{(Inkscape) Color is used for the text in Inkscape, but the package 'color.sty' is not loaded}%
    \renewcommand\color[2][]{}%
  }%
  \providecommand\transparent[1]{%
    \errmessage{(Inkscape) Transparency is used (non-zero) for the text in Inkscape, but the package 'transparent.sty' is not loaded}%
    \renewcommand\transparent[1]{}%
  }%
  \providecommand\rotatebox[2]{#2}%
  \ifx\svgwidth\undefined%
    \setlength{\unitlength}{240.80102608bp}%
    \ifx\svgscale\undefined%
      \relax%
    \else%
      \setlength{\unitlength}{\unitlength * \real{\svgscale}}%
    \fi%
  \else%
    \setlength{\unitlength}{\svgwidth}%
  \fi%
  \global\let\svgwidth\undefined%
  \global\let\svgscale\undefined%
  \makeatother%
  \begin{picture}(1,0.3487198)%
    \put(0,0){\includegraphics[width=\unitlength]{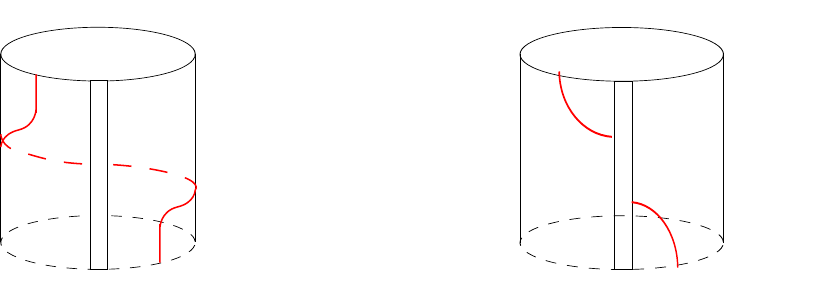}}%
    \put(0.13633449,0.17872568){\color[rgb]{0,0,0}\makebox(0,0)[lb]{\smash{$n$}}}%
    \put(0.71607562,0.00268422){\color[rgb]{0,0,0}\makebox(0,0)[lb]{\smash{$n$}}}%
    \put(0.76333355,0.15900511){\color[rgb]{0,0,0}\makebox(0,0)[lb]{\smash{$n+\varepsilon$}}}%
    \put(0.71607562,0.3283528){\color[rgb]{0,0,0}\makebox(0,0)[lb]{\smash{$n$}}}%
    \put(0.27101423,0.18549096){\color[rgb]{0,0,0}\makebox(0,0)[lb]{\smash{$=\underset{\varepsilon}{\sum}L_{\varepsilon}(n,r)$}}}%
  \end{picture}%
\endgroup%

\end{center}
 where $L_+(n,r)=(-1)^{n+1} e^{\frac{i \pi (n+2)}{r}} \left(\frac{\langle n+1 \rangle }{\langle n \rangle}\right)^{\frac{1}{2}}$
\newline and $L_{-}(n,r)=(-1)^{n+1} e^{\frac{-i \pi(n-2)}{r}} \left(\frac{\langle n-1 \rangle }{\langle n \rangle}\right)^{\frac{1}{2}}$
\newline
\newline All these coefficients are of the required form $\overline{c}(\gamma)F(\frac{c}{r} ,\frac{1}{r})$ 
\newline where $F$ is a smooth function defined on $V_{\gamma}= \lbrace (\tau ,\hbar ) \ / \ \tau_e \pm \hbar I_e \in U \rbrace$, and $\overline{c}(\gamma)$ is a sign factor, corresponding to taking the holonomy $\textrm{Hol}(\varepsilon , \gamma)$ of the cocyle $\overline{c}$ with values in $\lbrace \pm 1\rbrace$ along $\gamma$ already defined in Section \ref{sec:matcoef} and displayed in Figure \ref{fig:cocycle}.
\newline
\newline
\begin{figure}[htb]
  \begin{center}
  \def\svgwidth{190pt}
  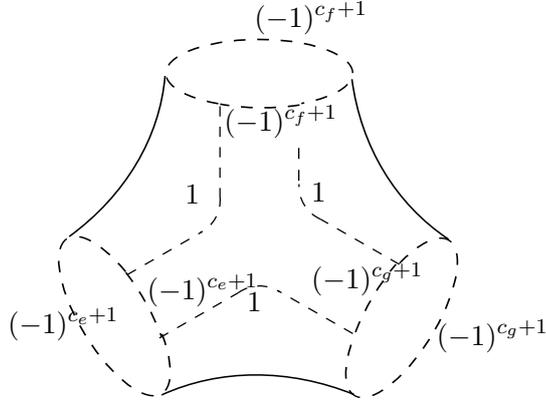
  \caption{The cocycle $\overline{c}$ on the pants bounded by the curves $C_e$, $C_f$ and $C_g$}
  \label{fig:cocycle}
  \end{center}
  \end{figure}
 If we have many curves in a pants or annulus, we only need to choose an order to make the fusions, and apply the latter formulas. For example, in the case of the pants , we obtain:
\newline
\newline 
\begingroup%
  \makeatletter%
  \providecommand\color[2][]{%
    \errmessage{(Inkscape) Color is used for the text in Inkscape, but the package 'color.sty' is not loaded}%
    \renewcommand\color[2][]{}%
  }%
  \providecommand\transparent[1]{%
    \errmessage{(Inkscape) Transparency is used (non-zero) for the text in Inkscape, but the package 'transparent.sty' is not loaded}%
    \renewcommand\transparent[1]{}%
  }%
  \providecommand\rotatebox[2]{#2}%
  \ifx\svgwidth\undefined%
    \setlength{\unitlength}{412.98981903bp}%
    \ifx\svgscale\undefined%
      \relax%
    \else%
      \setlength{\unitlength}{\unitlength * \real{\svgscale}}%
    \fi%
  \else%
    \setlength{\unitlength}{\svgwidth}%
  \fi%
  \global\let\svgwidth\undefined%
  \global\let\svgscale\undefined%
  \makeatother%
  \begin{picture}(1,0.30909845)%
    \put(0,0){\includegraphics[width=\unitlength]{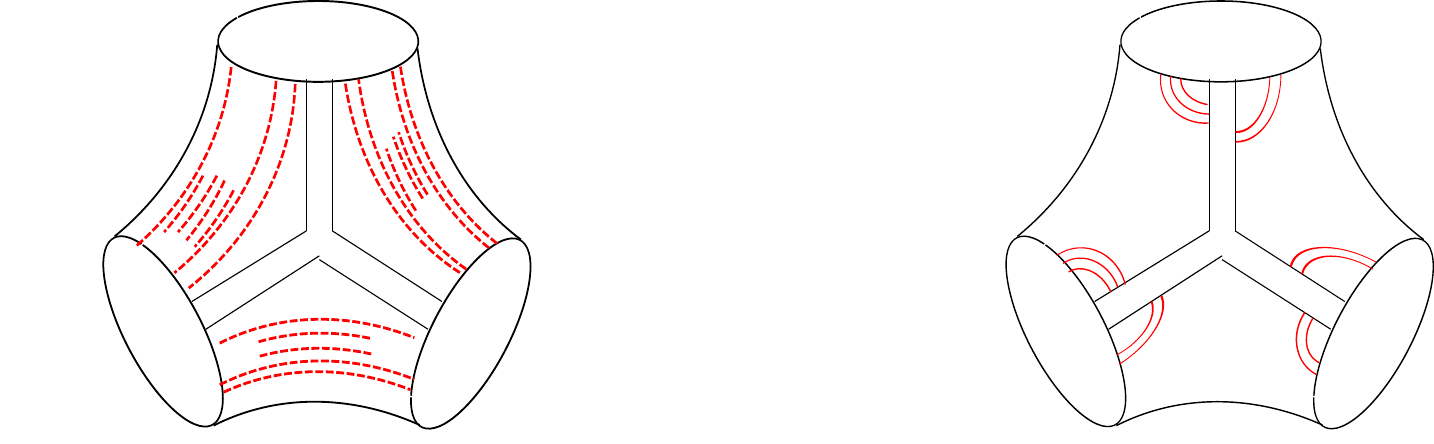}}%
    \put(0.20162273,0.26350342){\color[rgb]{0,0,0}\makebox(0,0)[lb]{\smash{$b$}}}%
    \put(0.30816288,0.069794){\color[rgb]{0,0,0}\makebox(0,0)[lb]{\smash{$a$}}}%
    \put(0.08539709,0.069794){\color[rgb]{0,0,0}\makebox(0,0)[lb]{\smash{$c$}}}%
    \put(0.82149274,0.27318886){\color[rgb]{0,0,0}\makebox(0,0)[lb]{\smash{$b$}}}%
    \put(0.9474038,0.07947949){\color[rgb]{0,0,0}\makebox(0,0)[lb]{\smash{$a$}}}%
    \put(0.72463803,0.069794){\color[rgb]{0,0,0}\makebox(0,0)[lb]{\smash{$c$}}}%
    \put(0.14350991,0.20539058){\color[rgb]{1,0,0}\makebox(0,0)[lb]{\smash{$\ldots$}}}%
    \put(-0.00069047,0.2217088){\color[rgb]{1,0,0}\makebox(0,0)[lb]{\smash{$\alpha$ curves}}}%
    \put(0.32804525,0.21681688){\color[rgb]{1,0,0}\makebox(0,0)[lb]{\smash{$\gamma$ curves}}}%
    \put(0.16288084,0.00199574){\color[rgb]{1,0,0}\makebox(0,0)[lb]{\smash{$\beta$ curves}}}%
    \put(0.76337992,0.19570511){\color[rgb]{0,0,0}\makebox(0,0)[lb]{\smash{$b+B$}}}%
    \put(0.88929098,0.1375923){\color[rgb]{0,0,0}\makebox(0,0)[lb]{\smash{$a+ A$}}}%
    \put(0.82149274,0.08916494){\color[rgb]{0,0,0}\makebox(0,0)[lb]{\smash{$c+C$}}}%
    \put(0.41502953,0.18349471){\color[rgb]{0,0,0}\makebox(0,0)[lb]{\smash{$=\underset{\varepsilon , \mu , \nu}{\sum}P_{\varepsilon , \mu , \nu}(a,b,c,r)$}}}%
  \end{picture}%
\endgroup%

\newline
\newline
where we wrote $A=\underset{i=1}{\overset{\beta +\gamma}{\sum}} \varepsilon_i$, $B = \underset{j=1}{\overset{\alpha+\gamma}{\sum}} \mu_j$, and $C=\underset{k=1}{\overset{\alpha+\beta}{\sum}} \nu_k$
\newline
\newline Here we have first used fusion on the $\alpha$ curves that go from $C_b$ to $C_c$, then the $\beta$ curves that run from $C_a$ to $C_c$, and finally the $\gamma$ curves from $C_a$ to $C_c$. With this order for the fusions, the coefficients $P_{\varepsilon ,\mu ,\nu}(a,b,c,r)$ are products of three factors corresponding to each serie of fusions: 
\newline 
\begin{scriptsize} 
\begin{eqnarray*}
F_{\mu_1 ,\nu_1}(a,b,c,r)F_{\mu_2 ,\nu_2}(a,b+\mu_1 ,c+\nu_1 ,r) \ldots F_{\mu_{\alpha} ,\nu_{\alpha}}(a, b+\sum_{i=1}^{\alpha-1} \mu_i ,c+\sum_{i=1}^{\alpha-1} \nu_i ,r)
\\ F_{\nu_{\alpha+1} ,\varepsilon_1}(b+\sum_{i=1}^{\alpha} \mu_i,a,c+\sum_{i=1}^{\alpha} \nu_i,r) \ldots F_{\nu_{\alpha+\beta} ,\varepsilon_{\beta}}(b+\sum_{i=1}^{\alpha} \mu_i ,a+\sum_{i=1}^{\beta-1}\varepsilon_i, c+\sum_{i=1}^{\alpha+\beta-1} \nu_i ,r)
\\  F_{\mu_{\alpha+1} ,\varepsilon_{\beta+1}}(c+\sum \nu, b+\sum_{i=1}^{\alpha} \mu_i,a+\sum_{i=1}^{\beta} \varepsilon_i ,r) \ldots F_{\mu_{\alpha+\gamma} ,\varepsilon_{\beta+\gamma}}(c+\sum \nu ,b+\sum_{i=1}^{\alpha+\gamma-1} \mu_i ,a+\sum_{i=1}^{\beta+\gamma-1} \varepsilon_i,r)
\end{eqnarray*}
\end{scriptsize}
\newline Notice that at every step of the fusion, the shifts in the color $c_e$ are sums of $\pm 1$ terms, one term for each arc intersecting $C_e$ that has been merged with $\Gamma$. Hence the coefficients $P_{\varepsilon ,\mu ,\nu}$ is defined and smooth on the required domain $V_{\gamma}= \lbrace (\tau ,\hbar) \ / \ \tau_e \pm I_e^{\gamma} \hbar \in U \rbrace$. Furthermore, at the end the shift in $c_e$ is no greater than the number of curves that intersect $C_e$ and of same parity as this number. 
\newline
\newline We now only need to explain what happens when we glue together two candlesticks.
\newline First, remark that we can only paste candlesticks with the same number of legs, and the same bottom color n. Moreover, if we paste two candlesticks $C(n , \varepsilon ,\Theta)$ and $C(n, \mu ,\Theta')$ with $\sum_{j} \mu_j \neq \sum_i \varepsilon_i$, then we obtain always 0 (as the vector space $V_r (\Sigma )$ of a sphere $\Sigma$ with two points marked by different colors is 0).
\newline
\newline
 \begin{prop}\label{gluing} The gluing of two candlesticks with $k$ legs $C(n,\varepsilon ,\Theta)$ et $C(n, \mu ,\Theta')$  with $\sum_{i=1}^{k} \varepsilon_i =\sum_{j=1}^{k} \mu_j$ is proportionnal to a band colored by $n+\sum \varepsilon_i$ joining the two boundary components of the annulus with no twist, the proportionality constant being $G(\frac{n}{r},\frac{1}{r})$, where G is a smooth function on $\lbrace (\tau ,\hbar) \ / \ \tau \pm  k \hbar \in (0,1) \rbrace$.
\end{prop} We should point out that in this proposition, the function G depends on $\Theta$, $\Theta'$, $\varepsilon$ and $\mu$.
\newline \textbf{Proof : } We prove this proposition by recurrence on the number of legs of the candlestick.
If we paste two candlesticks with only one leg, it is direct from the fusion rule that eliminate bigons (see Figure \ref{fig:fusion}), adding only a factor $(\frac{\langle c \pm 1 \rangle}{\langle c \rangle})^{\frac{1}{2}}$. Now, if $n=2$, the only delicate case is when the legs of the two part are disposed as in the third case of the figure below :
\newline
\newline
\begingroup%
  \makeatletter%
  \providecommand\color[2][]{%
    \errmessage{(Inkscape) Color is used for the text in Inkscape, but the package 'color.sty' is not loaded}%
    \renewcommand\color[2][]{}%
  }%
  \providecommand\transparent[1]{%
    \errmessage{(Inkscape) Transparency is used (non-zero) for the text in Inkscape, but the package 'transparent.sty' is not loaded}%
    \renewcommand\transparent[1]{}%
  }%
  \providecommand\rotatebox[2]{#2}%
  \ifx\svgwidth\undefined%
    \setlength{\unitlength}{338.18959968bp}%
    \ifx\svgscale\undefined%
      \relax%
    \else%
      \setlength{\unitlength}{\unitlength * \real{\svgscale}}%
    \fi%
  \else%
    \setlength{\unitlength}{\svgwidth}%
  \fi%
  \global\let\svgwidth\undefined%
  \global\let\svgscale\undefined%
  \makeatother%
  \begin{picture}(1,0.2805941)%
    \put(0,0){\includegraphics[width=\unitlength]{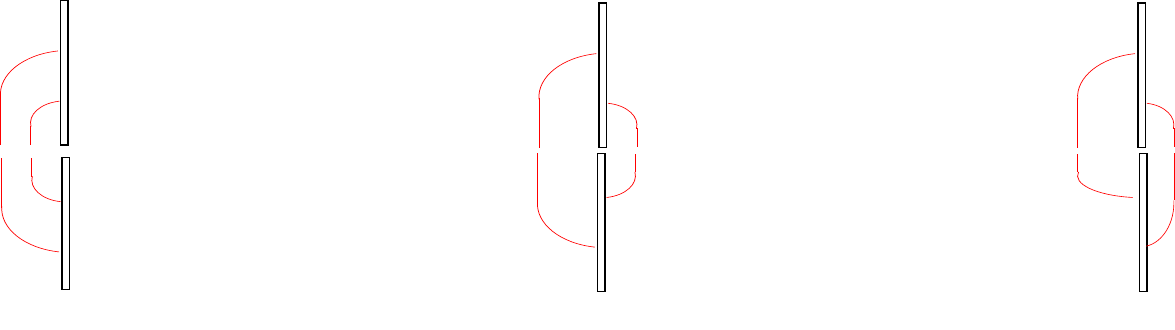}}%
    \put(0.02928854,0.00094573){\color[rgb]{0,0,0}\makebox(0,0)[lb]{\smash{Case 1}}}%
    \put(0.48969907,0.00016276){\color[rgb]{0,0,0}\makebox(0,0)[lb]{\smash{Case 2}}}%
    \put(0.95610676,0.00016276){\color[rgb]{0,0,0}\makebox(0,0)[lb]{\smash{Case 3}}}%
  \end{picture}%
\endgroup%

\newline
\newline
Indeed, in case 1 and 2, we could simply eliminate two bigons. For the case 3, we use the following switching legs formulas:
\newline
\newline
\begingroup%
  \makeatletter%
  \providecommand\color[2][]{%
    \errmessage{(Inkscape) Color is used for the text in Inkscape, but the package 'color.sty' is not loaded}%
    \renewcommand\color[2][]{}%
  }%
  \providecommand\transparent[1]{%
    \errmessage{(Inkscape) Transparency is used (non-zero) for the text in Inkscape, but the package 'transparent.sty' is not loaded}%
    \renewcommand\transparent[1]{}%
  }%
  \providecommand\rotatebox[2]{#2}%
  \ifx\svgwidth\undefined%
    \setlength{\unitlength}{449.74440603bp}%
    \ifx\svgscale\undefined%
      \relax%
    \else%
      \setlength{\unitlength}{\unitlength * \real{\svgscale}}%
    \fi%
  \else%
    \setlength{\unitlength}{\svgwidth}%
  \fi%
  \global\let\svgwidth\undefined%
  \global\let\svgscale\undefined%
  \makeatother%
  \begin{picture}(1,0.36066931)%
    \put(0,0){\includegraphics[width=\unitlength]{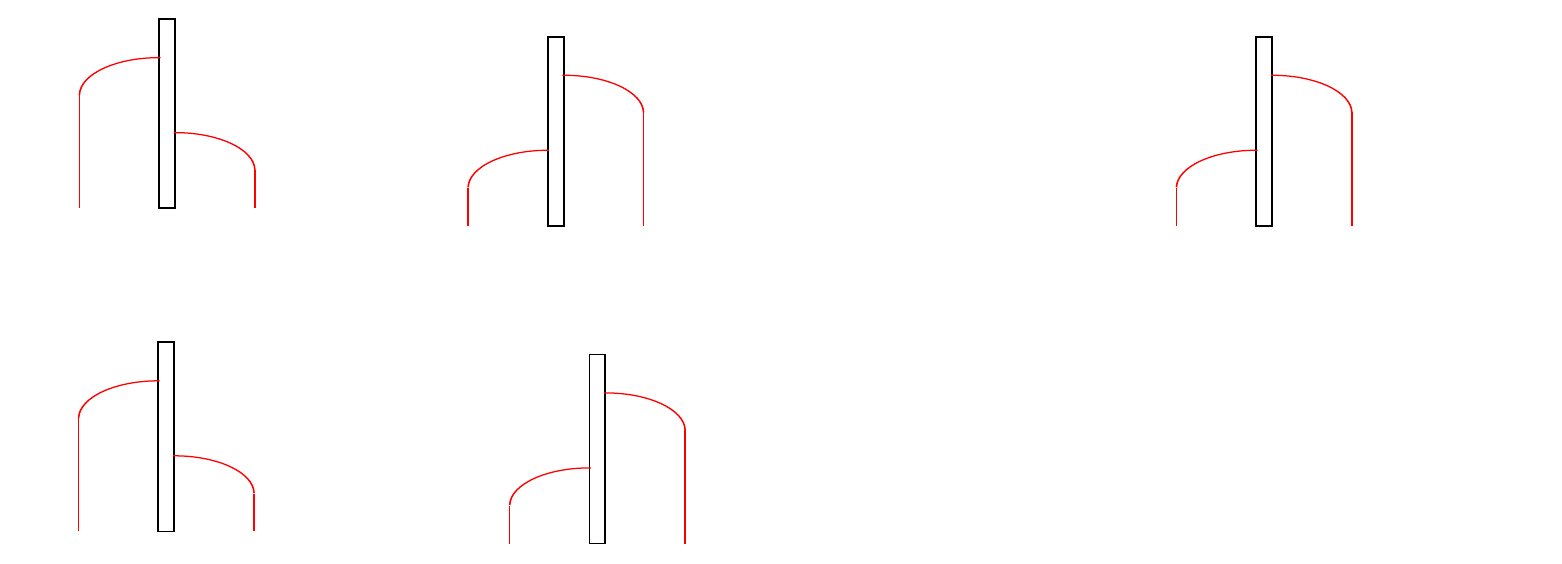}}%
    \put(0.17341879,0.24132413){\color[rgb]{0,0,0}\makebox(0,0)[lb]{\smash{$= \mp  \frac{\langle 1 \rangle}{\langle c \rangle}$}}}%
    \put(0.47581274,0.24132413){\color[rgb]{0,0,0}\makebox(0,0)[lb]{\smash{$+ \frac{(\langle c+1 \rangle \langle c-1 \rangle)^{1/2}}{\langle c \rangle}$}}}%
    \put(-0.00072822,0.12527385){\color[rgb]{0,0,0}\makebox(0,0)[lb]{\smash{$c \pm 2$}}}%
    \put(0.08099206,0.00421781){\color[rgb]{0,0,0}\makebox(0,0)[lb]{\smash{$c$}}}%
    \put(0.12185221,0.08492184){\color[rgb]{0,0,0}\makebox(0,0)[lb]{\smash{$c \pm 1$}}}%
    \put(0.20899453,0.04565747){\color[rgb]{0,0,0}\makebox(0,0)[lb]{\smash{$=$}}}%
    \put(0.31572184,0.11680897){\color[rgb]{0,0,0}\makebox(0,0)[lb]{\smash{$c \pm 2$}}}%
    \put(0.35702931,0.00166295){\color[rgb]{0,0,0}\makebox(0,0)[lb]{\smash{$c$}}}%
    \put(0.31572184,0.08123323){\color[rgb]{0,0,0}\makebox(0,0)[lb]{\smash{$c \pm 1$}}}%
    \put(0.12005514,0.34805142){\color[rgb]{0,0,0}\makebox(0,0)[lb]{\smash{$c $}}}%
    \put(0.0815121,0.21118334){\color[rgb]{0,0,0}\makebox(0,0)[lb]{\smash{$c$}}}%
    \put(0.12237224,0.29188737){\color[rgb]{0,0,0}\makebox(0,0)[lb]{\smash{$c \pm 1$}}}%
    \put(0.34240364,0.34805142){\color[rgb]{0,0,0}\makebox(0,0)[lb]{\smash{$c $}}}%
    \put(0.33039157,0.20507024){\color[rgb]{0,0,0}\makebox(0,0)[lb]{\smash{$c$}}}%
    \put(0.28903998,0.27689993){\color[rgb]{0,0,0}\makebox(0,0)[lb]{\smash{$c \pm 1$}}}%
    \put(0.79599459,0.34805142){\color[rgb]{0,0,0}\makebox(0,0)[lb]{\smash{$c $}}}%
    \put(0.78398253,0.20507019){\color[rgb]{0,0,0}\makebox(0,0)[lb]{\smash{$c$}}}%
    \put(0.74263099,0.27689993){\color[rgb]{0,0,0}\makebox(0,0)[lb]{\smash{$c \mp 1$}}}%
  \end{picture}%
\endgroup%

\begin{figure}
  \begin{center}
  \def\svgwidth{330pt}
  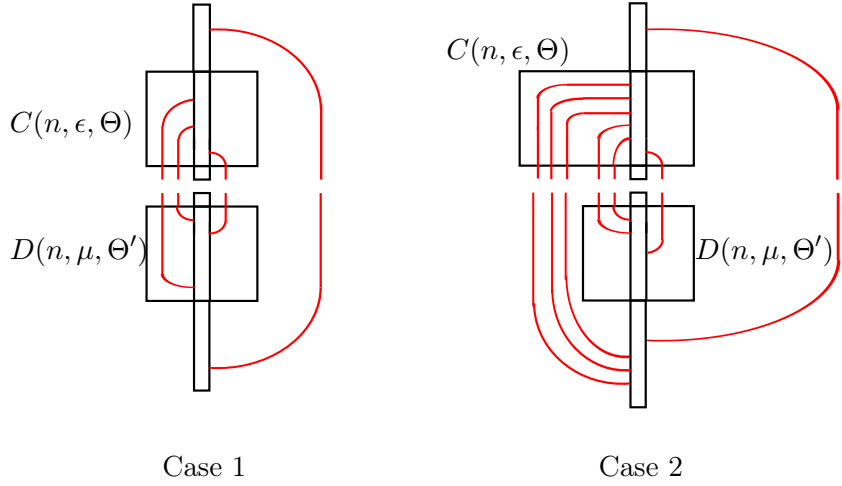
  \caption{The two cases of pasting candlesticks with $k$ legs}
  \label{fig:pastchand}
    \end{center}
  \end{figure}

\newline
\newline
\newline This shows Proposition \ref{gluing} for $k \leq 2$.
\newline 
\newline Now, suppose we glue two candlesticks with $k+1$ legs. We have two cases as in Figure \ref{fig:pastchand}:
\newline
\newline
In case 1, the upper leg of the upper candlestick and the bottom leg of the bottom candlestick both go to the right (or both to the left), the gluing is obtained by gluing two candlesticks with k legs, then suppressing a bigon. The factor we get is of the form $G(\frac{n}{r},\frac{1}{r}) \left(\frac{\langle n + \sum_{i=1}^{k+1} \varepsilon_i \rangle}{\langle n + \sum_{i=1}^{k} \varepsilon_i \rangle}\right)^{\frac{1}{2}}$, which is indeed a function of $(\frac{n}{r},\frac{1}{r})$ that is smooth on the domain we claimed.
\newline On the contrary, in case 2, the upper leg of the upper part and the bottom leg of the bottom part go to different sides, but by applying the formulas to switch legs, we can reduce this to the former case.
\subsection{Examples of $\psi$-symbol}
\label{sec:ex}
\begin{figure}
  \centering
  \def\svgwidth{200pt}
  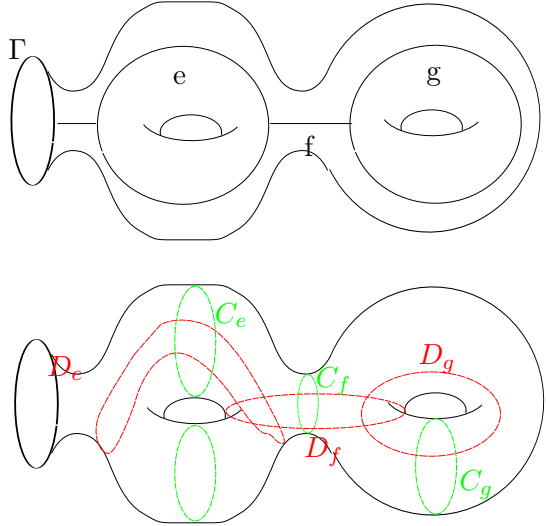
  \caption{The family of curves $D_e$ on $\Sigma$ associated to the trivalent banded graph $\Gamma$}
  \label{fig:pantalon}
\end{figure}

We derive expressions of $\psi$-symbol for two families of curves on
$\Sigma$: the first family consists of the curves $C_e$ of the pants
decompostion itself, and the other of curves $(D_e )_{e \in E}$ that are
in some sense dual to the curves $C_e$. The $D_e$ are defined this way: if
$e$ is an internal edge that joints a vertex to itself, then $D_e$ is a loop
parallel to $e$. If $e$ joints two different vertices, then $D_e$ consists
of two arcs parallel to $e$ that we close into a loop as in Figure \ref{fig:pantalon}.
\newline 
Note that $C_e$ and $D_f$ intersect each other if and only if
$e=f$, and in this case they intersect once or twice, and that $p(C_e)$ and $p(D_e)$ are null homotopic in $H_1(\Gamma , \partial \Gamma ,\mathbb{Z}/2)$.
\newline
\newline 
\begin{prop}\label{calcul} We have for any $e \in E$ and $c \in U_r$:
\newline 1) $T_r^{C_e} \varphi_c = -2 \cos (\pi \frac{c_e}{r}) \varphi_c$ et
$\sigma^{C_e}=-2 \cos (\pi \frac{c_e}{r})[0]$
\newline 2) In the case where $e$ is an edge joining an internal vertex to itself as in the picture:
\newline
\newline
\begin{center}
\begingroup%
  \makeatletter%
  \providecommand\color[2][]{%
    \errmessage{(Inkscape) Color is used for the text in Inkscape, but the package 'color.sty' is not loaded}%
    \renewcommand\color[2][]{}%
  }%
  \providecommand\transparent[1]{%
    \errmessage{(Inkscape) Transparency is used (non-zero) for the text in Inkscape, but the package 'transparent.sty' is not loaded}%
    \renewcommand\transparent[1]{}%
  }%
  \providecommand\rotatebox[2]{#2}%
  \ifx\svgwidth\undefined%
    \setlength{\unitlength}{191.30118354bp}%
    \ifx\svgscale\undefined%
      \relax%
    \else%
      \setlength{\unitlength}{\unitlength * \real{\svgscale}}%
    \fi%
  \else%
    \setlength{\unitlength}{\svgwidth}%
  \fi%
  \global\let\svgwidth\undefined%
  \global\let\svgscale\undefined%
  \makeatother%
  \begin{picture}(1,0.5798876)%
    \put(0,0){\includegraphics[width=\unitlength]{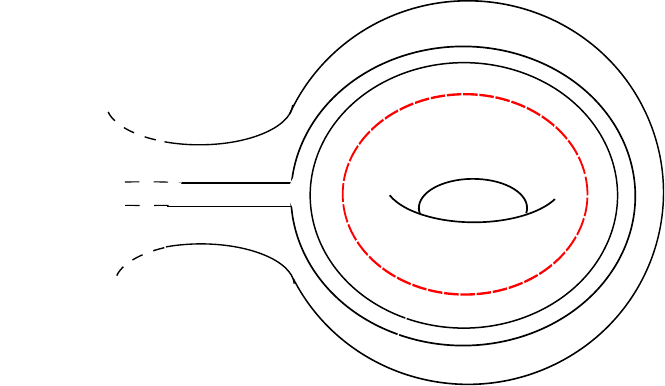}}%
    \put(0.77131663,0.51650855){\color[rgb]{0,0,0}\makebox(0,0)[lb]{\smash{$c_e$}}}%
    \put(0.24579537,0.23624501){\color[rgb]{0,0,0}\makebox(0,0)[lb]{\smash{$c_f$}}}%
  \end{picture}%
\endgroup%

\end{center}
 we have $\sigma^{D_e}(\frac{c}{r} ,\theta ,\frac{1}{r})= \left(W(\pi \frac{c_e}{r}, \pi \frac{c_f}{r}, \frac{\pi}{r}) e^{i \theta_e} +W(\pi \frac{c_e}{r} ,\pi \frac{c_f}{r}, - \frac{\pi}{r}) e^{-i \theta_e}\right)[0]$ 
\newline where $W(\tau ,\alpha ,\hbar) =\left( \frac{\sin(\tau +\alpha/2 +\hbar/2) \sin(\tau -\alpha/2+\hbar/2) }{\sin(\tau)\sin(\tau+\hbar)} \right)^{\frac{1}{2}}$
\newline 3) In the case where $e$ is an edge between two distinct internal vertices as in the picture:
\newline
\newline
\begin{center}
\begingroup%
  \makeatletter%
  \providecommand\color[2][]{%
    \errmessage{(Inkscape) Color is used for the text in Inkscape, but the package 'color.sty' is not loaded}%
    \renewcommand\color[2][]{}%
  }%
  \providecommand\transparent[1]{%
    \errmessage{(Inkscape) Transparency is used (non-zero) for the text in Inkscape, but the package 'transparent.sty' is not loaded}%
    \renewcommand\transparent[1]{}%
  }%
  \providecommand\rotatebox[2]{#2}%
  \ifx\svgwidth\undefined%
    \setlength{\unitlength}{170.475bp}%
    \ifx\svgscale\undefined%
      \relax%
    \else%
      \setlength{\unitlength}{\unitlength * \real{\svgscale}}%
    \fi%
  \else%
    \setlength{\unitlength}{\svgwidth}%
  \fi%
  \global\let\svgwidth\undefined%
  \global\let\svgscale\undefined%
  \makeatother%
  \begin{picture}(1,0.74175099)%
    \put(0,0){\includegraphics[width=\unitlength]{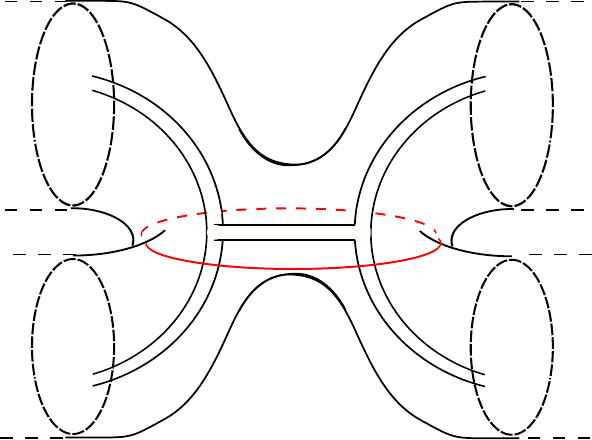}}%
    \put(0.2345893,0.63136275){\color[rgb]{0,0,0}\makebox(0,0)[lb]{\smash{$c_a$}}}%
    \put(0.21876725,0.18923978){\color[rgb]{0,0,0}\makebox(0,0)[lb]{\smash{$c_d$}}}%
    \put(0.44620966,0.42827372){\color[rgb]{0,0,0}\makebox(0,0)[lb]{\smash{$c_e$}}}%
    \put(0.71518503,0.65472693){\color[rgb]{0,0,0}\makebox(0,0)[lb]{\smash{$c_b$}}}%
    \put(0.71122958,0.22338755){\color[rgb]{0,0,0}\makebox(0,0)[lb]{\smash{$c_c$}}}%
    \put(0.43408124,0.21314125){\color[rgb]{1,0,0}\makebox(0,0)[lb]{\smash{$D_e$}}}%
  \end{picture}%
\endgroup%

\end{center}
we have $\sigma^{D_e} (\frac{c}{r} ,\theta , \frac{1}{r}) = -\left(I(\pi \tau ,\pi \hbar) +J(\pi \tau , \pi \hbar) e^{2 i \theta_e} +J(\pi (\tau -2 \hbar \delta_e) ,\pi  \hbar) e^{-2 i \theta_e}\right)[0]$ 
\newline Here, we have set $\tau =\frac{c}{r}$, $\hbar =\frac{1}{r}$, $\delta_e$ for the Kronecker symbol,  
\begin{eqnarray*}
I(\tau ,\hbar) &=& 2 \cos (\tau_c +\tau_d -\hbar)
\\ & & +4\frac{\sin(\frac{\tau_a +\tau_d -\tau_e -\hbar}{2})\sin(\frac{\tau_a -\tau_d +\tau_e +\hbar}{2})\sin(\frac{\tau_b +\tau_c -\tau_e -\hbar}{2})\sin(\frac{\tau_b -\tau_c +\tau_e +\hbar}{2}) }{\sin(\tau_e)\sin(\tau_e +\hbar)}
\\ & & +4\frac{\sin(\frac{\tau_a +\tau_d +\tau_e -\hbar}{2})\sin(\frac{-\tau_a +\tau_d +\tau_e -\hbar}{2})\sin(\frac{\tau_b +\tau_c +\tau_e -\hbar}{2})\sin(\frac{-\tau_b +\tau_c +\tau_e -\hbar}{2}) }{\sin(\tau_e)\sin(\tau_e -\hbar)}
\end{eqnarray*}
 and 
 \begin{eqnarray*}
 J(\tau , \hbar) &=& 4 (\frac{\sin(\frac{\tau_a +\tau_d -\tau_e -\hbar}{2})\sin(\frac{\tau_a -\tau_d +\tau_e +\hbar}{2})\sin(\frac{\tau_b +\tau_c -\tau_e -\hbar}{2})\sin(\frac{\tau_b -\tau_c +\tau_e +\hbar}{2}) }{\sin(\tau_e)\sin(\tau_e +\hbar)}
 \\ & & \frac{\sin(\frac{\tau_a +\tau_d +\tau_e +\hbar}{2})\sin(\frac{-\tau_a +\tau_d +\tau_e +\hbar}{2})\sin(\frac{\tau_b +\tau_c +\tau_e +\hbar}{2})\sin(\frac{-\tau_b +\tau_c +\tau_e +\hbar}{2}) }{\sin(\tau_e +\hbar)\sin(\tau_e +2\hbar)})
\end{eqnarray*}
\end{prop}

The expressions of $T_r^{C_e}$ and $T_r^{D_e}$ can be derived by using fusion rules. The computations are rather long in the last case, but straighforward.
\newline  These expressions, as well as the expressions of the $\psi$-symbol of the curves $C_e$ and $D_e$ were already given in \cite{MP}. They also checked by hand that the formulas of Theorem \ref{maintheo2} were satisfied by these curves. We will only derive from the formulas that the zero-th and first order term for these curves are related as in Theorem \ref{maintheo2}, as this is rather quick, and will come of use later:
\newline
\newline 
\begin{prop}\label{calc2} Let $\gamma$ be any of the curves $C_e$ or $D_e$. 
\newline Then $\sigma^{\gamma}(\tau ,\theta ,\hbar)=\sigma^{\gamma}(\tau ,\theta,0) +\frac{\hbar}{2 i}\underset{e \in E}{\sum} \frac{\partial^{2}}{\partial \tau_e \partial \theta_e} \sigma^{\gamma}(\tau ,\theta ,0)+o(\hbar)$
\end{prop}
 \textbf{Proof :} For $C_e$, there is not much to prove: as $\sigma^{C_e}$ does not depend on $\hbar$, the first order term vanishes, and as $\sigma^{C_e}$ does not depend on $\theta_e$, $\frac{\partial^{2}}{\partial \tau_e \partial \theta_e} \sigma^{\gamma}(\tau ,\theta ,0)$ also vanishes.
\newline For the curves $D_e$, we need to separate the case 1 where $e$ joints a vertex to itself, and the case 2 where it joints two distinct vertices.
\newline In case 1 depicted by the figure above, we have 
\newline $\sigma^{D_e}(\tau ,\theta ,\hbar)= (W(\pi \tau_e, \pi \tau_f, \pi \hbar) e^{i \theta_e} +W(\pi \tau_e ,\pi \tau_f, - \pi \hbar) e^{-i \theta_e})[0]$. Notice that from the formula of $W$ given above, we get $W(\pi \tau_e ,\pi \tau_f ,\pi \hbar)=W(\pi (\tau_e+\frac{\hbar}{2}),\pi \tau_f ,0)+o(\hbar)$. Thus
\begin{eqnarray*}\sigma^{D_e}(\tau,\theta,\hbar) &=& \sigma^{D_e}(\tau,\theta,0)+\frac{\hbar}{2}(W(\pi \tau_e,\pi \tau_f,0)e^{i \theta_e} -W(\pi \tau_e,\pi \tau_f,0)e^{-i \theta_e})[0]+o(\hbar)
\\ &=& \sigma^{D_e}(\tau ,\theta,0) +\frac{\hbar}{2 i}\underset{e \in E}{\sum} \frac{\partial^{2}}{\partial \tau_e \partial \theta_e} \sigma^{D_e}(\tau ,\theta ,0)+o(\hbar)
\end{eqnarray*}
as expected.
\newline Finally, in the case 2 above, we have 
\newline $\sigma^{D_e} (\tau ,\theta , \hbar) = -(I(\pi \tau ,\pi \hbar ) +J(\pi \tau , \pi \hbar) e^{2 i \theta_e} +J(\pi (\tau -2 \hbar \delta_e) , \pi \hbar) e^{-2 i \theta_e})[0]$. It is easily seen that $J(\tau ,\hbar)=J(\tau +\hbar \delta_e ,0)$ for $\delta_e$ Kronecker symbol. Thus the only thing we need is to prove that $I(\tau ,\hbar)=I(\tau ,0)+o(\hbar)$. This is a bit more tricky:
\newline First, notice that we can write 
\begin{displaymath}I(\tau,\hbar)=2 \cos(\tau_c+\tau_d-\hbar) +\frac{1}{\sin (\tau_e)} (F(\tau_e+\hbar)-F(-\tau_e+\hbar)) +o(\hbar)
\end{displaymath}
 where
\begin{eqnarray*}F(\tau_e) &=& 4\frac{\sin(\frac{\tau_a+\tau_d-\tau_e}{2})\sin(\frac{\tau_a-\tau_d+\tau_e}{2})\sin(\frac{\tau_b+\tau_c-\tau_e}{2})\sin(\frac{\tau_b-\tau_c+\tau_e}{2})}{\sin(\tau_e)}
\\ &=& \frac{(\cos(\tau_d-\tau_e)-\cos(\tau_a))(\cos(\tau_c-\tau_e)-\cos(\tau_b))}{\sin(\tau_e)}
\end{eqnarray*}
Therefore, the first order term of $I(\tau,\hbar)$ is  $\hbar \left(2\sin(\tau_c+\tau_d) +\frac{2}{\sin (\tau_e)} \frac{d}{d \tau_e} \mathcal{P}(F)(\tau_e)\right)$, where $\mathcal{P}(F)$ is the even part of the function $F$. From the formula above, we have
\begin{displaymath}\mathcal{P}(F)(\tau_e)=\sin(\tau_c+\tau_d)\cos(\tau_e) -\cos(\tau_a)\sin(\tau_c)-\cos(\tau_b)\sin(\tau_d)
\end{displaymath}
So that $\frac{1}{\sin (\tau_e)} \frac{d}{d \tau_e} \mathcal{P}(F)(\tau_e)=-\sin(\tau_c+\tau_d)$, and the first order of $I(\tau,\hbar)$ vanishes.
\newline
\newline The computations of $\sigma^{C_e}$ and $\sigma^{D_e}$ were previously used by \cite{MP} to prove a version of Theorem \ref{maintheo2} for the punctured torus and the 4-holed sphere. Their approach was to derive from the formulas  that the asymptotic estimate of Theorem \ref{maintheo2} was valid for the curves $C_e$, $D_e$ and $\tau_{C_e}(D_e)$ where $\tau_{C_e}$ denotes the Dehn-twist along $C_e$. Then they used the compatibility of the $\psi$-symbol with the product in $K(\Sigma ,-e^{\frac{i \pi \hbar}{2}})$ to prove that if Theorem \ref{maintheo2} is verified by $\gamma$ and $\delta $ two multicurves, then it is also true for their product $\gamma \cdot \delta$. This yielded Theorem \ref{maintheo2} for all multicurves in the punctured torus and the 4-holed sphere, as the curves $C_e$, $D_e$, and $\tau_{C_e}(D_e)$ were sufficient to generate the Kauffman algebra.
\newline
However, this approach fails in higher genus, as this set of curves no longer generate the Kauffman algebra. Therefore, we developped an other approach to tackle the higher genus cases, which was also more conceptual and required less computations. Our fundamental idea is to use the multiplicativity of the $\psi$-symbol together with the theorem of Bullock (reminded in Section \ref{sec:matcoef}) to view the zero-th and first order term of the $\psi$-symbol in terms of algebra morphism or derivation of algebras on $\textrm{Reg}(\mathcal{M}'(\Sigma))$. We then only need to compare this general shape with the values of the $\psi$-symbol on a few curves to get the formula of Theorem \ref{maintheo2}. (In fact, we will only need the values on the $C_e$ for the zero order term, while the first order term also required the values on $D_e$)
\section{Principal symbol and representation spaces}
\label{sec:prsym}
This section will be centered on the study of the principal symbol $\sigma^{\gamma}(\tau , \theta , 0)$, that is the zero-order of the $\psi$-symbol $\sigma^{\gamma}(\tau , \theta , \hbar)$. The goal of this paragraph is to establish the formula for the principal symbol, which is stated our main theorem: $\sigma_{\chi}^{\gamma}(\tau , \theta ,0)=f_{\gamma}(R_{\chi}(\tau ,\theta ))$, where $f_{\gamma}$ is the function on $\mathcal{M}(\Sigma)$ such that $f_{\gamma}(\rho)=-\textrm{Tr}(\rho(\gamma))$, and $R_{\chi}$ are action-angles parametrization on $\mathcal{M}(\Sigma)$.
\subsection{The intersection algebra $A_{\Gamma}$ and multiplicativity of the $\psi$-symbol}
\label{sec:multi}
The aim of this section is to prove the property of compatibility of the $\psi$-symbol with the multiplication in $K(\Sigma ,-1)$, given by Proposition \ref{mult}:
\newline
\newline
\textbf{Proposition \ref{mult} : }Let $\gamma$ and $\delta$ be two multicurves on $\Sigma$. Then we have the following asymptotic expression:
\begin{displaymath} \sigma^{\gamma \cdot \delta}(\tau , \theta , \hbar ) =\left(\sigma^{\gamma}(\tau , \theta , \hbar) \sigma^{\delta}(\tau , \theta , \hbar ) + \frac{\hbar}{i} \underset{e}{\sum} \partial_{\tau_{e}}\sigma^{\gamma}(\tau , \theta , \hbar ) \partial_{\theta_{e}} \sigma^{\delta}(\tau , \theta , \hbar ) \right) + o(\hbar) 
\end{displaymath}
 A version of this proposition appeared in \cite{MP}, but they worked with an other definition of the $\psi$-symbol, which took values in $\mathbb{C}$, whereas in our definition, the $\psi$-symbol takes values in $A_{\Gamma}$.
\newline Thus, before we turn to the proof of the proposition, we would like to start this section by pointing out a few facts about the algebras $A_{\Gamma}$.
\newline
\newline We know that $A_{\Gamma}$ is a commutative $\mathbb{C}$-algebra of dimension $2^{g}$, and that implies that it is isomorphic to the algebra $\mathbb{C}^{2^{g}}$. We would like to give a more explicit form to the isomorphism.
\newline First, we point out that $A_{\Gamma}$ can be viewed as a quotient of a group algebra. Indeed, take the set $G=H_1(\Gamma ,\partial \Gamma ,\mathbb{Z}/2) \times \lbrace \pm 1 \rbrace$, and define the product on $G$ by:
\begin{displaymath} (\alpha , \varepsilon) \cdot (\beta , \mu) = (\alpha +\beta , \varepsilon \mu \langle \alpha , \beta \rangle)
\end{displaymath}
Here, $\langle \cdot , \cdot \rangle$ is the intersection form in $H_1 (\Gamma , \partial \Gamma , \mathbb{Z}/2)$. Then  $A_{\Gamma}$ is clearly the quotient $\mathbb{C}[G]/(\alpha,-1) \sim - (\alpha , 1)$. The algebra $\mathbb{C}[G]$ is isomorphic to the algebra $\mathbb{C}^{\hat{G}}$, where $\hat{G}=\textrm{Hom}(G ,\lbrace \pm 1 \rbrace)$ and the components of the isomorphism are of the form
\begin{displaymath} \underset{g \in G}{\sum} \alpha_g g \rightarrow  \underset{g \in  G}{\sum} \alpha_g \rho (g)
\end{displaymath}
where $\rho$ goes over $\hat{G}$. Then the components of the isomorphism $A_{\Gamma} \rightarrow \mathbb{C}^{2^{g}}$ are given by 
\begin{displaymath} \underset{\gamma \in H_1 (\Gamma ,\partial \Gamma ,\mathbb{Z}/2)}{\sum} \alpha_{\gamma} [\gamma] \rightarrow  \underset{\gamma \in  H_1(\Gamma ,\partial \Gamma ,\mathbb{Z}/2)}{\sum} \alpha_{\gamma} \rho (\gamma , 1)
\end{displaymath}
where $\rho$ goes over the set $\hat{A}_{\Gamma}$ of representations $G \rightarrow \lbrace \pm 1 \rbrace$ such that $\rho (0 ,-1) =-1$. We point out that the quotients of two such representations corresponds to representations $H_1(\Gamma ,\partial \Gamma ,\mathbb{Z}/2) \rightarrow \lbrace \pm 1 \rbrace$, and $\hat{A}_{\Gamma}$ has a structure of affine space over $H_{1}(\Gamma ,\partial \Gamma ,\mathbb{Z}/2)$. 
\newline
\newline Another way of thinking of this affine space structure is to view $\hat{A}_{\Gamma}$ as the set of "relative spin-structures" over $(\Gamma,\partial \Gamma)$. Indeed, define a linear map $\chi :A_{\Gamma} \mapsto \mathbb{C}$ by $\chi([\gamma])=(-1)^{q(\gamma)}$. Then we have $\chi \in \hat{A}_{\Gamma}$ if and only if for any $\gamma,\delta \in H_{1}(\Gamma,\partial \Gamma ,\mathbb{Z}/2)$,we have $q(\gamma+\delta)=q(\gamma)+q(\delta)+\langle \gamma, \delta \rangle$, where $\langle \cdot ,\cdot \rangle$ is the intersection form in $H_1(\Gamma,\partial \Gamma, \mathbb{Z}/2)$. 
\newline
\newline Now the principal symbol $\sigma^{\gamma}(\tau ,\theta , 0)$ is a morphism $K(\Sigma ,-1) \rightarrow A_{\Gamma}$, we also label his components $\sigma_{\chi}^{\gamma}(\tau ,\theta) =\chi(\sigma^{\gamma}(\tau ,\theta ,0))$ for every $\chi \in \hat{A}_{\Gamma}$.
\newline
\newline \textbf{Proof of Proposition \ref{mult}: }We fix $r > 0$ and we take $\gamma$ and $\delta$ two multicurves on $\Sigma$. According to Theorem \ref{maintheo1}, the matrix coefficents of the operator $T_{r}^{\gamma}$ can be written as:
\begin{displaymath} T_{r}^{\gamma} \varphi_{c}=\overline{c}(\gamma) \underset{k :E \rightarrow \mathbb{Z}}{\sum} F_k^{\gamma}(\tau ,\hbar) \varphi_{c+k}
\end{displaymath}
with the $F_k^{\gamma}$ being smooth functions on $V_{\gamma}$ such that $F_{k}^{\gamma}=0$ as soon as there is some $e \in E$ such that $|k_e|>I_e^{\gamma}$ or $k_e \not\equiv I_e^{\gamma} (\textrm{mod} \ 2)$.
\newline As $\gamma \in K(\Sigma ,-e^{\frac{i \pi }{2 r}}) \rightarrow T_r^{\gamma}\in \textrm{End}(V_r (\Sigma))$ is an morphism of algebras, we have:
\begin{displaymath}T_r^{\gamma \cdot \delta}\varphi_c = T_r^{\gamma} (T_r^{\delta} \varphi_c )
\end{displaymath} 
writing $\tau =\frac{c}{r}$ and $\hbar=\frac{1}{r}$, and from the above expression of the matrix coefficients, we get 
\begin{eqnarray*}T_r^{\gamma \cdot \delta} \varphi_c &=& \underset{m : E \rightarrow \mathbb{Z}}{\sum}\left( \underset{k+l=m}{\sum}F_l^{\gamma}(\tau +k \hbar ,\hbar)F_k^{\delta}(\tau ,\hbar) \overline{c}(\gamma) \overline{c+k}(\gamma) \right) \varphi_{c+m}
\\ &=& \overline{c}(\gamma) \overline{c}(\delta) i(\gamma ,\delta) \underset{m :E \rightarrow \mathbb{Z}}{\sum}\left( \underset{k+l=m}{\sum}F_l^{\gamma}(\tau +k \hbar ,\hbar)F_k^{\delta}(\tau ,\hbar) \right) \varphi_{c+m}
\end{eqnarray*}
Here, we could factor $\overline{c+k}(\gamma)$ out of the sum, as if $F_k^{\delta} \neq 0$ then $k_e \equiv I_e^{\delta} [2]$ and $\overline{k}(\gamma)=i(\gamma ,\delta)=\underset{e \in E}{\prod}(-1)^{I_e^{\delta}(C_e^{' *}(\gamma)+C_e^{'' *}(\gamma))}$. We later show that $i(\gamma ,\delta)$ correspond to the intersection in $H_1 (\Gamma ,\partial \Gamma , \mathbb{Z}/2)$ of the projections of $\gamma$ and $\delta$.
\newline Now as $K(\Sigma ,-e^{\frac{i \pi \hbar}{2}})$ is generated by multicurves, we can write $\gamma \cdot \delta =\underset{\lambda}{\sum} f_{\lambda}(\hbar )\lambda$, and in this sum, $f_{\lambda} \neq 0$ only when $[\lambda] =[\gamma]+[\delta] \in H_{1}(\Sigma , \mathbb{Z}/2)$, according to the Kauffman relations. When it is the case, we have $\overline{c}(\lambda) =\overline{c}(\gamma) \overline{c}(\delta)$. We can write another formula for the curve operator of the product:
\begin{displaymath}T_r^{\gamma \cdot \delta} \varphi_c =\underset{m}{\sum}\left(\underset{\lambda}{\sum} \overline{c}(\lambda)f_{\lambda}(\hbar) F_{m}^{\lambda}(\tau ,\hbar) \right) \varphi_{c+m}	
\end{displaymath}
So, identifying coefficients in the two formulas, we get:
\begin{displaymath}\underset{\lambda}{\sum} f_{\lambda}(\hbar) F_m^{\lambda}(\tau ,\hbar) =\left(\underset{k+l=m}{\sum}F_l^{\gamma}(\tau +k \hbar ,\hbar)F_k^{\delta}(\tau ,\hbar)\right) i(\gamma ,\delta)
\end{displaymath}
Now, remember that we defined the $\psi$-symbol of an arbitrary element of $K(\Sigma ,-e^{\frac{i \pi \hbar}{2}})$ by extending lineary the formula for multicurves. Thus, we have:
\begin{displaymath}\sigma^{\gamma \cdot \delta} (\tau , \theta ,\hbar)=\underset{m}{\sum}\underset{\lambda}{\sum}f_{\lambda}(\hbar) F_m^{\lambda}(\tau ,\hbar) e^{i m \theta}[\lambda]
\end{displaymath} 
recalling that $[\lambda] =[\gamma]+[\delta]$ and using the previous identity of coefficients:
\begin{eqnarray*} &=& i(\gamma ,\delta)\underset{m}{\sum}\left(\underset{k+l=m}{\sum}F_l^{\gamma}(\tau +k \hbar ,\hbar)F_k^{\delta}(\tau ,\hbar)\right) e^{i m \theta} [\gamma +\delta]
\\ &=& i(\gamma , \delta) \langle p_*(\gamma) , p_*(\delta)\rangle \left(\sigma^{\gamma}(\tau ,\theta ,\hbar)\sigma^{\delta}(\tau ,\theta ,\hbar)+\frac{\hbar}{i}\underset{e \in E}{\sum}\partial_{\tau_e}\sigma^{\gamma} (\tau ,\theta ,\hbar) \partial_{\theta_e} \sigma^{\delta}(\tau ,\theta ,\hbar)\right) +o(\hbar)
\end{eqnarray*}
All that is left to prove is to show that the formula given for $i(\gamma ,\delta)$ actually computes the intersection in $H_{1}(\Gamma ,\partial \Gamma ,\mathbb{Z}/2)$.
\newline Let $\gamma$ and $\delta$ be two curves in $\Gamma$. We can isotope $\gamma$ and $\delta$ so that $\delta$ lies in the interior of $\Gamma$, and $\gamma$ follows the edges of the cell decomposition of $\Gamma$. Then the intersection points lie only in the curves $p(C_e)=L_e$. The number of intersection points of $\gamma$ and $\delta$ in $L_e$ is congruent modulo 2 to $\sharp (\delta \cap L_e) L_e^{*}(\gamma)$ where $L_e^{*}$ is the dual to the cell $L_e$. If $\tilde{\gamma}$ and $\tilde{\delta}$ are lifts of $\gamma$ and $\delta$ to $\Sigma$, then $L_e^{*}(\gamma)=C_e^{' *}(\tilde{\gamma})+C_e^{''*}(\tilde{\gamma})$ and $\sharp (\delta \cap L_e) =\sharp (\tilde{\delta} \cap C_e)$ modulo 2, hence the formula for $i(\tilde{\gamma},\tilde{\delta})$ computes the number of intersection points of $\gamma$ and $\delta$ modulo 2.
\subsection{Principal symbol and the $\textrm{SL}_2$-character variety}
\label{sec:charac}
This section aims to elaborate a link between the components of the principal symbol $\sigma_{\chi}$ and functions on the space of representations $\pi_1 (\Sigma) \rightarrow \textrm{SL}_{2}(\mathbb{C})$.
\newline We will start our study of the principal symbol by the following proposition, which describes which values $\sigma_{\chi}^{\gamma}(\tau , \theta)$ can take:  
\newline 
\newline
\begin{prop}\label{major} For any multicurve $\gamma$ and $\chi \in \hat{A}_{\Gamma}$, we have:
\newline 1)  $\sigma_{\chi}^{\gamma}(\tau , \theta) \in \mathbb{R}$
\newline 2)  $|\sigma_{\chi}^{\gamma}(\tau , \theta)| \leq 2^{n(\gamma)}$ where $n(\gamma)$ is the number of components of $\gamma$.
\end{prop}
\textbf{Proof : } 1) We recall that the components of the $\psi$-symbol $\sigma_{\chi}^{\gamma}$ are complex-valued. The stated property comes from the fact that curve operators are Hermitian: for any multicurve $\gamma$, and every $r$, the operator $T_{r}^{\gamma}$ is a Hermitian endomorphism of $V_{r}(\Sigma)$. By definition, we have
$T_{r}^{\gamma} \varphi_{c} =\underset{k}{\sum} F_{k}^{\gamma}(\frac{c}{r} , \frac{1}{r}) \varphi_{c+k}$. As the basis $(\varphi_{c})_{c \in U_{r}}$ is a Hermitian basis, we get $F_{-k}^{\gamma}(\frac{c+k}{r} , \frac{1}{r})=\overline{F_{k}^{\gamma}(\frac{c}{r} , \frac{1}{r})}$ for all $c \in U_{r}$. Then for $r \rightarrow +\infty$ we have $F_{-k}^{\gamma}(\tau , 0)=\overline{F_{k}^{\gamma}(\tau , 0)}$. 
\newline Hence $\sigma_{\chi}^{\gamma}(\tau , \theta)=\chi (\gamma) \underset{k}{\sum}F_{k}^{\gamma}(\tau , 0) e^{i k \cdot \theta}  \in \mathbb{R}$, $\forall (\tau ,\theta ) \in U \times (\mathbb{R}/2 \pi \mathbb{Z})^{E}$
\newline
\newline 2) We want to find a majoration of $|\sigma_{\chi}^{\gamma}(\tau , \theta)|$, where $\gamma$ is a multicurve. By definition, we have $\sigma_{\chi}^{\gamma}(\tau , \theta)=\chi (\gamma)\underset{k}{\sum}F_{k}^{\gamma}(\tau , 0 ) e^{i k \cdot \theta}$. From one hand, we know that the coefficients $F_{k}^{\gamma}$ are zero as soon as there is $e$ such that $|k_{e}| > I_{e}^{\gamma}=\sharp (\gamma \cap C_{e})$. The number of non-zero coefficients is then lower than $M_{\gamma}=\underset{e \in E}{\prod} (2 I_{e}^{\gamma} +1)$. On the other hand, for any $r \geq 2$ and $c \in U_{r}$:
\begin{displaymath} F_{k}^{\gamma}(\frac{c}{r},\frac{1}{r})=\langle T_{r}^{\gamma} \varphi_{c} , \varphi_{c+k} \rangle \leq ||T_{r}^{\gamma}||
\end{displaymath}
 We reminded in the preliminary section the spectral radius of $T_r^{\gamma}$ is always $\leq 2^{n(\gamma)}$. Thus we have $|F_{k}^{\gamma}(\frac{c}{r} , \frac{1}{r} ) | \leq 2^{n(\gamma)}$ for every r > 0 and every $c \in U_{r}$. Taking the limit, we get $|F_{k}^{\gamma}(\tau , 0)| \leq 2^{n(\gamma)}$.
\newline These two estimations only allow us to write $|\sigma_{\chi}^{\gamma}(\tau , \theta)| \leq M_{\gamma} 2^{n(\gamma)}$. To obtained to the promised inequality, we use the multiplicativity of $\sigma_{\chi}^{\cdot}(\tau , \theta)$: 
\newline We have, for any integer $p$: $|\sigma_{\chi}^{\gamma^{p}}(\tau , \theta)| = |\sigma_{\chi}^{\gamma}(\tau , \theta )|^{p}$.
But $\gamma^{p}$ is also a multicurve, obtained by taking $p$ parallel copies of each component of $\gamma$.
\newline
We deduce that $|\sigma_{\chi}^{\gamma^{p}}(\tau , \theta)| \leq M_{\gamma^{p}} 2^{n(\gamma^{p})} \leq p^{3 g -3} M_{\gamma} 2^{p n(\gamma)}$ and so we get by taking $p \rightarrow +\infty$ that $|\sigma_{\chi}^{\gamma}(\tau , \theta )| \leq 2^{n(\gamma)}$ for all 
$(\tau , \theta ) \in U \times (\mathbb{R}/2 \pi \mathbb{Z})^{E}$.
\newline
\newline
Now, recall that the components of the $\psi$-symbol
\begin{displaymath} \sigma_{\chi}(\tau , \theta) \ : \ K(\Sigma ,-1) \rightarrow \mathbb{C}
\end{displaymath}
are morphisms of algebras. There is a simple description of all such morphism of algebras: indeed, by the Theorem \ref{regfunc} that we recalled in Section \ref{sec:matcoef}, we have an isomorphism
 \begin{displaymath}K(\Sigma , -1) \simeq \textrm{Reg}(\mathcal{M}'(\Sigma))
 \end{displaymath} 
 where $\mathcal{M}'(\Sigma)$ stands for $\textrm{Hom}(\pi_{1} \Sigma , SL_{2}(\mathbb{C}))//SL_{2}(\mathbb{C})$, the spaces of caracters of the fundamental group of $\Sigma$ in $SL_{2}(\mathbb{C})$. This space is an affine algebraic variety, and we are writing $\textrm{Reg}(\mathcal{M}'(\Sigma))$ for the set of regular functions from $\mathcal{M}'(\Sigma)$ to $\mathbb{C})$. A morphism of algebras $\phi$ from $\textrm{Reg}(\mathcal{M}'(\Sigma))$ to $\mathbb{C}$ is always of the form 
\begin{displaymath} \phi \ : \ f \in \textrm{Reg}(\mathcal{M}'(\Sigma))  \ \rightarrow \ f(\rho)
\end{displaymath}
for some $\rho \in \mathcal{M}'(\Sigma)$. We deduce the existence of maps 
\begin{displaymath} R_{\chi} \ : \ U \times (\mathbb{R}/2 \pi \mathbb{Z})^{E} \rightarrow \mathcal{M}'(\Sigma)
\end{displaymath}
such that $\sigma_{\chi}^{\gamma}(\tau , \theta)=f_{\gamma}(R_{\chi}(\tau , \theta ))$.
\subsection{A system of actions-angles coordinates on the $\textrm{SU}_{2}$-character variety} 
\label{sec:action}
This paragraph will be devoted to study the maps $R_{\chi}$ more closely, the aim is to prove that it actually gives actions-angles coordinates on the character variety $\textrm{Hom}(\pi_1 (\Sigma) ,\textrm{SU}_{2})/\textrm{SU}_{2}$ that we will note $\mathcal{M}(\Sigma)$. 
\newline In $\mathcal{M}(\Sigma)$ there is an open dense subset $\mathcal{M}_{\textrm{irr}}(\Sigma)$ consisting of all conjugacy of irreducible representations. It is a well-known fact that $\mathcal{M}_{\textrm{irr}}(\Sigma)$ has a structure of (smooth) symplectic variety therefore we call it the smooth part of $\mathcal{M}(\Sigma)$.
\newline The maps $R_{\chi}$ are at first sight their image in $\mathcal{M}'(\Sigma)$. Again, we have a subset $\mathcal{M}^{'}_{\textrm{irr}}(\Sigma) \subset \mathcal{M}'(\Sigma)$ consisting of conjugacy classes of irreducible representations, and there is a structure of (smooth) symplectic variety on it (that restricts to that of $\mathcal{M}_{\textrm{irr}}(\Sigma)$.   
\newline We have two remarks:
\newline First, we point out that $R_{\chi}$ is always an non-commutative representation. Indeed, for a commutative representation, we would have for $e$,$f$,$g$ three adjacent edges
\begin{displaymath}h_{C_e}(\rho)+h_{C_f}(\rho)=h_{C_g}(\rho)
\end{displaymath}
for one of the three ordering of $e$,$f$,$g$, or have $h_{C_e}(\rho)+h_{C_f}(\rho)+h_{C_g}(\rho)=2$. This can not happen for $R_{\chi}(\tau ,\theta)$ as $(h_{C_e})_{e \in E}$ maps it to $\tau \in U$, and we have strict inequalities $\tau_g <\tau_e+\tau_f$, $\tau_e+\tau_f+\tau_g<2$.
\newline 
\newline 
 Our second point is that the map $R_{\chi}$ is smooth. This has a sense as by our first remark  its image is in the smooth part of $\mathcal{M}'(\Sigma)$. as stated in Definition 2.1, $(\tau , \theta ) \rightarrow \sigma^{\gamma}( \tau , \theta , 0)$ is a smooth function on $U \times (\mathbb{R}/2 \pi \mathbb{Z})^{E}$, for all $\gamma \in K(\Sigma , -1)$. So $(\tau ,\theta) \rightarrow \textrm{Tr}(R_{\chi}(\tau , \theta) (\gamma))$ is smooth for every $\gamma \in \pi_1(\Sigma)$. As the space $\mathcal{M}'(\Sigma)$ can be parametrized by a finite collection of coordinates $\rho \rightarrow \textrm{Tr}(\rho(\gamma_j))$, where $\gamma_j \in \pi_1(\Sigma)$, the map $R_{\chi} : U \times (\mathbb{R}/2 \pi \mathbb{Z})^{E} \rightarrow \mathcal{M}'(\Sigma)$ is smooth. 
\newline
\newline
\begin{prop}\label{values}The maps $R_{\chi}$ take in fact values in $\mathcal{M}_{\textrm{irr}}(\Sigma)=\textrm{Hom}(\pi_{1} \Sigma , \textrm{SU}_{2})/\textrm{SU}_{2}$ 
\end{prop}
\textbf{Proof : }Indeed, we have seen with Proposition \ref{major} that $\sigma_{\chi}^{\gamma}(\tau ,\theta )$ is real-valued. We can use a well-known lemma:
\newline 
\newline 
\textbf{Lemma : } Any irreductible subgroup $G \subset \textrm{SL}_{2}(\mathbb{C})$ such that the trace of all elements of G are real is conjugated to either a subgroup of $\textrm{SL}_{2}(\mathbb{R})$ or a subgroup of $\textrm{SU}_{2}$. 
\newline 
\newline The proof of this lemma is based only on elementary algebra, manipulating trace of products of elements of G. A detailed proof can be found for example in \cite[p.3040-3041]{HK}. 
\newline As we have $\sigma^{\gamma}(\tau , \theta ,0)=- \textrm{Tr}(R(\tau ,\theta)(\gamma)) \in \mathbb{R}$, we get that $R(\tau ,\theta )$ is conjugated to either a representation in $\textrm{SL}_{2}(\mathbb{R})$ or a representation in $\textrm{SU}_{2}$.
\newline 
\newline To prove Proposition \ref{values}, we still need to dismiss the case where the image of $R_{\chi}(\tau , \theta)$ would be conjugated to a subgroup of $\textrm{SL}_{2}(\mathbb{R})$. To this end, we use the point 2) of Proposition \ref{major}, who states that 
$|\textrm{Tr}(R_{\chi}(\tau ,\theta)\gamma)| \leq 2$ for every $\gamma \in \pi_{1}(\Sigma)$ representing a simple closed curve on $\Sigma$. We use the following lemma, proved in \cite{GKM}:
\newline
\newline
\textbf{Lemma 2.1 : }Let $\rho \ : \ \pi_{1}(\Sigma) \rightarrow \textrm{PSL}_{2}(\mathbb{C})$  be a non-elementary representation, then there exists two simple loops $a$ and $b$ intersecting once such that $\rho (a)$ and $\rho (b)$ are loxodromic ($|\textrm{Tr}(\rho (a))|>2$, $|\textrm{Tr}(\rho (b))|>2$) and non commuting.
\newline
\newline
This lemma appeared as Lemma 3.1.1 in \cite{GKM} and follows from elementary considerations in hyperbolic geometry. From the lemma, we get that since $R(\tau ,\theta)(a)$ is never loxodromic, it must be an elementary representation into $\textrm{PSL}_2(\mathbb{C})$. But if $R(\tau ,\theta)$ was conjugated to a representation in $\textrm{SL}_2(\mathbb{R})$, it would be a commutative representation, and we saw that $R(\tau, \theta)$ is not.
\newline
\newline
\begin{prop}\label{act-ang} For any $\chi \in \hat{A}_{\Gamma}$, the map 
\newline $R_{\chi} \ : \ (\tau , \theta ) \in U \times (\mathbb{R}/2 \pi \mathbb{Z})^{E} \rightarrow R_{\chi}(\tau ,\theta) \in \mathcal{M}(\Sigma)$ gives action-angle coordinates on the symplectic variety $\mathcal{M}_{\textrm{irr}}(\Sigma)$.
\end{prop}
We remind that when a pants decomposition $\mathcal{C}=\lbrace C_{e} \rbrace_{e \in E}$ of $\Sigma$ is given, the family of functions $h_{C_{e}}=\frac{1}{\pi} \textrm{Acos}(-\frac{f_{C_{e}}}{2})$ constitutes a moment mapping $h \ : \ \mathcal{M}(\Sigma) \rightarrow \overline{U}$. The variables $\tau_{e}$ are the action coordinates associated to this moment mapping:
\begin{displaymath}h_{C_{e}}(R_{\chi}(\tau , \theta))= \frac{1}{\pi} \textrm{Acos}(-\frac{f_{C_{e}}(R_{\chi}(\tau , \theta))}{2})=\frac{1}{\pi} \textrm{Acos}(-\frac{\sigma_{\chi}^{C_{e}}(\tau , \theta)}{2})= \tau_{e}
\end{displaymath}
where the third equality comes from the computation of the operator $T_{r}^{C_{e}}$ given in Section \ref{sec:comput}: for any coloration $c$ of $E$, we have $T_{r}^{C_{e}} \varphi_{c} = - 2 \cos (\frac{\pi c}{r}) \varphi_{c}$, so that $\sigma_{\chi}^{C_{e}}(\tau , \theta , \hbar)=F_{0}^{C_{e}}(\tau , \hbar )\chi([0])= - 2 \cos( \pi \tau_{e}) $. 
\newline
\newline The fact that $(\tau , \theta )$ is a system of action-angle coordinates on $\mathcal{M}(\Sigma)$ can be described the following way:
\begin{displaymath} \omega = \underset{e \in E}{\sum} d \tau_{e} \wedge d \theta_{e}
\end{displaymath}
where $\omega$ refers to the symplectic form on the variety $\mathcal{M}(\Sigma)$.
\newline It also amounts to the fact that the vector fields $\partial_{\theta_{e}}$ and $X_{h_{C_{e}}}$ (the symplectic gradient associated to the function $h_{C_{e}}$) on $\mathcal{M}(\Sigma)$ are equals. This equality of vector fields can be rewritten in terms of Poisson brackets:
\begin{displaymath} \forall f \in C^{\infty}(\mathcal{M}(\Sigma) ,\mathbb{C}),\ \forall \tau , \theta \ \textrm{we have } \ \lbrace h_{C_{e}} , f \rbrace = \frac{\partial}{\partial \theta_{e}} f(R_{\chi}(\tau,\theta))
\end{displaymath}
We only need to verify this equality when $f$ is one of the function $f_{\gamma}$, where $\gamma \in K(\Sigma ,-1)$ as the Poisson bracket is a first order differential operator, and any function f on $\mathcal{M}(\Sigma)$ can be approximated at order 1 by a trace function $f_{\gamma}$. By linearity, we can show it only for $\gamma$ a multicurve.
\newline To compute such Poisson brackets, we can apply the formula of Goldman \cite{G86} that we recalled in Section \ref{sec:matcoef}:
\newline We note $\varepsilon$ the linear map
\begin{displaymath}\varepsilon \ : \ K(\Sigma , -e^{\frac{i \pi \hbar}{2}}) \rightarrow K(\Sigma , -1) \simeq \textrm{Reg}(\mathcal{M}'(\Sigma))
\end{displaymath}
\begin{displaymath}
 \ \ \ \ \underset{\gamma \ \textrm{multicurve}}{\sum} c_{\gamma}(\hbar) \gamma \rightarrow \underset{\gamma \ \textrm{multicurve}}{\sum} c_{\gamma}(0) \gamma
\end{displaymath}  
For $\gamma$ and $\delta \in K(\Sigma , -e^{\frac{i \pi \hbar}{2}})$ we have:
\begin{displaymath} \lbrace f_{\varepsilon(\gamma)} , f_{\varepsilon(\delta)}\rbrace = f_{\varepsilon(\frac{i}{\hbar}[\gamma , \delta ])} 
\end{displaymath}
with $[\gamma , \delta ] = \gamma \cdot \delta - \delta \cdot \gamma \in K(\Sigma , -e^{\frac{i \pi \hbar}{2}})$.
\newline  We apply the above formula
to compute $\lbrace h_{C_{e}} , f_{\gamma} \rbrace$ for any $\gamma \in K(\Sigma , -e^{\frac{i \pi \hbar}{2}})$: We recall that $h_{C_{e}}=\frac{1}{\pi} \textrm{Acos}(-\frac{f_{C_{e}}}{2})$. Our strategy to compute the Poisson bracket is to approximate $h_{C_e}$ with polynomials in $f_{C_e}$.
\newline On a neighboorhood $V$ of $R_{\chi}(\tau ,\theta)$, $f_{\gamma}$ has values in an open set $(-2+\eta ,2-\eta) \subset [-2 ,2]$. We choose a sequence of polynomials $P_j$ such that $P_j$ converge to the map $x \mapsto \frac{1}{\pi}\textrm{Acos}(-\frac{x}{2})$ on $(-2+\eta ,2-\eta)$ in the  $C^{1}$-topology. The Poisson bracket being a differential operator of order one, we have that $\lbrace P_j (f_{C_e}),f_{\gamma}\rbrace$ converges uniformly on $V$ to $\lbrace h_{C_e} ,f_{\gamma} \rbrace$ when $j \rightarrow +\infty$.
\newline 
\newline 
Now, the maps $\lbrace \cdot , f_{\gamma} \rbrace \ : \ C^{\infty}(\mathcal{M}(\Sigma)) \rightarrow C^{\infty}(\mathcal{M}(\Sigma))$
\newline and $\frac{i}{\hbar}[\cdot , \gamma] \ : \ K(\Sigma ,-1) \rightarrow K(\Sigma ,-1)$ being derivations of algebras, we have by Goldman's formula: 
\begin{displaymath}\lbrace P_j (f_{C_{e}}) , f_{\gamma} \rbrace (R_{\chi}(\tau , \theta))=f_{\varepsilon(\frac{i}{\hbar}[P_j(C_{e}) , \gamma ])}(R_{\chi}(\tau , \theta))=
\sigma_{\chi}^{\varepsilon(\frac{i}{\hbar}[P_j(C_{e}) , \gamma ])}(\tau ,\theta , 0)
\end{displaymath}
We compute this last quantity: we recall that we wrote $T_{r}^{\gamma} \varphi_{c}= \underset{k}{\sum} F_{k}^{\gamma}(\tau , \hbar ) \varphi_{c+k}$ and we gave in Section \ref{sec:ex} the expression $T_{r}^{C_{e}} \varphi_{c} =-2 \cos (\pi \tau_{e}) \varphi_{c}$. 
Hence $T_{r}^{P_j(C_{e})} \varphi_{c} = P_j(-2 \cos(\pi \tau_{e})) \varphi_{c}$.
We deduce that for $c \in U_{r}$:
\begin{displaymath}T_{r}^{[P_j(C_{e}) ,\gamma]} \varphi_{c}= 
\underset{k}{\sum} P_j(-2\cos(\pi (\tau_e + k_e \hbar)))F_{k}^{\gamma}(\tau ,\hbar )\varphi_{c+k} -  
\underset{k}{\sum} P_j(-2\cos(\pi \tau_e)) F_{k}^{\gamma}(\tau ,\hbar )\varphi_{c+k} 
\end{displaymath}
so that as $[C_e^{k}]=[0]$ in $A_{\Gamma}$,
\begin{displaymath}\sigma_{\chi}^{\varepsilon(\frac{i}{\hbar}[P_j(C_e),\gamma])}(\tau,\theta,0)=\frac{i}{2\pi}\underset{k}{\sum}\left.\frac{P_j(-2 \cos(\pi (\tau_e +k_e \hbar)))-P_j(-2 \cos(\pi \tau_e))}{\hbar}\right|_{\hbar=0} F_k^{\gamma}(\tau,0) e^{i k\cdot \theta} \chi (\gamma)
\end{displaymath}
When $j$ tends to $+\infty$, as $P_j$ approach the function $x \mapsto \frac{1}{\pi}\textrm{Acos}(-\frac{x}{2})$, this quantity goes to
 \begin{displaymath}
\underset{k}{\sum} i k_e F_{k}^{\gamma}(\tau , \hbar) e^{i k \cdot \theta}\chi(\gamma)
=\frac{\partial}{\partial_{\theta_{e}}}\sigma_{\chi}^{\gamma}(\tau , \theta , 0)
=\frac{\partial}{\partial_{\theta_{e}}}f_{\gamma}(R_{\chi}(\tau , \theta ))
\end{displaymath}
The last equality ends the proof: we have indeed $\lbrace h_{C_{e}} , f_{\gamma} \rbrace (R_{\chi}(\tau ,\theta)) =  \frac{\partial}{\partial \theta_{e}} f_{\gamma}(R_{\chi}(\tau,\theta))$ for every multicurve $\gamma$ , and $R_{\chi}$ give a  action-angle parametrization of $\mathcal{M}_{\textrm{irr}}(\Sigma)$.
\newline
\newline 
\textbf{Origin of angle coordinates}
\newline
\newline Finally, we want to investigate how exactly $R_{\chi}$ varies with $\chi \in \hat{A}_{\Gamma}$. We recall that according to Section \ref{sec:multi}, the values of two different morphisms $\chi$ and $\chi '$ on $[\gamma]$ differ by a representation $\rho : H_{1}(\Gamma ,\partial \Gamma ,\mathbb{Z}/2) \rightarrow \lbrace \pm 1 \rbrace$.
\newline Let us also get more precise about angle coordinates. We recall that we have an hamiltonian $h : \mathcal{M}_{\textrm{irr}}(\Sigma) \rightarrow U$, given by $(h(\rho))_e=\frac{1}{\pi}\textrm{Acos}(-\frac{\textrm{Tr}(\rho(C_e))}{2})$. The hamiltonian flow gives an action of $\mathbb{R}^E$ on $\mathcal{M}_{\textrm{irr}}(\Sigma)$. This action has a kernel
\begin{displaymath} \Lambda=\textrm{Vect}_{\mathbb{Z}}\lbrace(2 \pi u_e )_{ e \in E}, \ \ \pi (u_e+u_f+u_g)_{(e,f,g) \in S}\rbrace
\end{displaymath}
where $(u_e)_{e \in E}$ is the canonical basis of $\mathbb{R}^E$, $E$ is the set of edges of $\Gamma$, and $S$ is the set of triple of edges adjacent to the same vertex in $\Gamma$. We also define $\Lambda '=\textrm{Vect}_{\mathbb{Z}}(\pi u_e) \supset \Lambda$. The quotient $\Lambda '/\Lambda$ then acts on $\mathcal{M}^{\textrm{irr}}(\Sigma)$ by $\pi u_e \cdot \rho (\gamma)= (-1)^{( C_e ,\gamma )} \rho(\gamma)$, where $( \cdot ,\cdot )$ is the intersection form in $\Sigma$.
\newline  
\newline Now that we know that the maps $R_{\chi}$ give action-angle coordinates on $\mathcal{M}_{\textrm{irr}}(\Sigma)$, the only ambiguity is the choice of the origin of the angle part. That is we must have for any $\chi ,\chi ' \in \hat{A}_{\Gamma}$ that $R_{\chi '}(\tau ,\theta)=R_{\chi}(\tau ,\theta +v_{\chi ,\chi'})$ for a fixed vector $v_{\chi ,\chi'} \in \mathbb{R}/\Lambda$.
\newline We use the values of $R_{\chi}$ on the curves $D_e$ to get the origin of angle coordinates. We have
\newline $\textrm{Tr}(R_{\chi}(\tau,\theta)(D_e))=-\sigma_{\chi}^{D_e}(\tau,\theta,0)=-2 W(\pi \tau ,0) \cos(\theta_e)$ if $e$ joins a vertex to itself, 
\newline or $=I(\pi \tau , 0) +2 J(\pi \tau ,0)\cos (2 \theta_e)$ otherwise. We see that in the first case, $\theta_e =0$ is the unique minimum of $\textrm{Tr}(R_{\chi}(\tau ,\theta)(D_e))$, so that this origin of this coordinate is the same for all $\chi \in \hat{A}_{\Gamma}$. In the second case, $\theta_e \mapsto \textrm{Tr}(R_{\chi}(\tau ,\theta)(D_e))$ has exactly two maxima, one for $\theta_e=0$, one for $\theta_e=\pi$. So $\theta$ is fixed modulo $\pi u_e$. Thus for $\chi,\chi' \in \hat{A}_{\Gamma}$, we have $v_{\chi ,\chi'} \in \Lambda ' /\Lambda$.
\newline Taking two elements $\chi ,\chi'$ in $\hat{A}_{\Gamma}$ we know that they differ by a morphism $\rho :H_1(\Gamma, \partial \Gamma ,\mathbb{Z}/2) \rightarrow \lbrace \pm 1 \rbrace$. It is possible to recover the vector $v_{\chi ,\chi'} \in \Lambda '/\Lambda$ from the representation $\rho$: indeed, by Poincare duality, one can write $\rho(p_*(\gamma))=(-1)^{\langle C ,\gamma \rangle }$ where $C \in H_1(\Sigma, \mathbb{Z}/2)$, $p_*$ is the projection $H_1(\Sigma ,\mathbb{Z}/2) \rightarrow H_1(\Gamma, \partial \Gamma ,\mathbb{Z}/2)$ and $\langle \cdot, \cdot \rangle$ is the intersection form in $H_1(\Sigma,\mathbb{Z}/2)$. Remember that $p_*$ maps each $C_e$ to zero, so that the intersection of $C$ with each $C_e$ must vanish. As the $C_e$ generate a Lagragian of $H_1(\Sigma,\mathbb{Z}/2)$, $C$ is a linear combination of the $C_e$ and this yields a vector $v_{\rho} \in \Lambda'/\Lambda$ such that $R_{\rho \chi}(\tau ,\theta)=R_{\chi}(\tau, \theta +v_{\rho})$.
\newline
\newline We need to note that when $\Gamma$ is a planar graph we can drop these complicated consideration of angle origins and we could have taken the $\psi$-symbol to be just $\mathbb{C}$-valued. Indeed, in this case the intersection form in $H_1(\Gamma, \partial \Gamma ,\mathbb{Z}/2)$ is trivial, and the image of $H_1(\Sigma,\mathbb{Z}/2) \rightarrow H_1(\Gamma, \partial \Gamma ,\mathbb{Z}/2)$ is $\lbrace 0 \rbrace$, so that all components of the $\psi$-symbol are the same.
\section{First order of the $\psi$-symbol}
\label{sec:frstor}
In this section, we investigate the first order term in $\hbar$ of the asymptotic expansion of the $\psi$-symbol. We identify this term by linking it with the principal symbol, of which we already know a formula. 
\newline We remind that for $\gamma$ a multicurve, the map $(\tau , \hbar , \theta) \rightarrow \sigma^{\gamma}(\tau , \theta , \hbar )$ is defined as a finite sum of smooth functions on $V_{\gamma} $, and $V_{\gamma}$ is a neighborhood of $U \times \lbrace 0 \rbrace$ in $U \times [0,1]$. We may write, for any multicurve $\gamma$:
\begin{displaymath}\sigma^{\gamma}(\tau , \theta , \hbar)=\sigma^{\gamma}(\tau , \theta ,0) +\hbar (\Delta_{\gamma}(\tau , \theta)+D_{\gamma}(\tau , \theta)) + o (\hbar)
\end{displaymath}
Here, $\Delta_{\gamma}(\tau ,\theta)$ refers to the expected first order as in Theorem 2: 
\newline $\Delta_{\gamma}(\tau ,\theta)=\frac{1}{2 i} \underset{e}{\sum} \frac{\partial^{2}}{\partial \tau_{e} \partial \theta_{e}} \sigma^{\gamma}(\tau , \theta ,0)$. Hence what we want to prove in this section is that the default $D_{\gamma}(\tau , \theta )$ is zero for all $\gamma$ and $(\tau , \theta) \in U \times (\mathbb{R}/2 \pi \mathbb{Z})^{E}$. 
\newline We remark that the previous expressions defines $\Delta(\tau ,\theta)$ and $D(\tau ,\theta)$ as maps from the set of multicurves to $A_{\Gamma}$, what we can extend by linearity to linear maps $K(\Sigma , -e^{\frac{i \pi \hbar}{2}}) \rightarrow A_{\Gamma}[[\hbar]]$. 
\newline Furthermore, $\Delta_{\gamma}$ and $D_{\gamma}$ are some linear combinations of partial derivatives of the smooth functions $F_{k}$ on $V_{\gamma}$, they are both smooth on $U \times (\mathbb{Z}/2 \pi \mathbb{Z})^{E}$.
\newline
\newline
\begin{prop}\label{default}
For any multicurve $\gamma$ and for all $(\tau ,\theta)$, the default $D_{\gamma}(\tau , \theta)$ vanishes, so that the first order term of $\sigma^{\gamma}(\tau , \theta ,\hbar )$ is
$\Delta_{\gamma}(\tau ,\theta)=\frac{1}{2 i} \underset{e}{\sum} \frac{\partial^{2}}{\partial \tau_{e} \partial \theta_{e}} \sigma^{\gamma}(\tau , \theta ,0)$. 
\end{prop}
The proof relies on the two following lemmas:
\newline
\newline
\begin{lem}\label{deriv}
Let $(\tau ,\theta)$ be in $U \times (\mathbb{R}/2 \pi \mathbb{Z})^{E}$. We will provide $\mathbb{C}$ with a structure of $K(\Sigma , -1)$-module (or equivalently of  $\textrm{Reg}(\mathcal{M}'(\Sigma))$-module): for $x \in \mathbb{C}$ and $f \in \textrm{Reg}(\mathcal{M}'(\Sigma))$, we define $f \cdot x = f(R_{\chi}(\tau , \theta)) x$. Then the corresponding component of the default $\gamma \mapsto \chi (D_{\gamma}(\tau ,\theta ))$ is a derivation of $K(\Sigma ,-1)$-modules from $K(\Sigma , -1)$ to $\mathbb{C}$.
\end{lem}
\begin{lem}\label{tang} For the same structure of $\textrm{Reg}(\mathcal{M}'(\Sigma))$-module on $\mathbb{C}$ as above, we have an isomorphism $\textrm{Der}(\textrm{Reg}(\mathcal{M}'(\Sigma)),\mathbb{C}) \simeq T_{R_{\chi}(\tau ,\theta )}\mathcal{M}(\Sigma)$ sending a vector $X \in T_{R_{\chi}(\tau ,\theta )}\mathcal{M}(\Sigma)$ to the derivation $f \rightarrow \mathcal{L}_{X} f (R_{\chi}(\tau ,\theta))$, and the vector fields $(\partial \tau_e ,\partial \theta_e)$ give basis of the tangent spaces $T_{R_{\chi}(\tau ,\theta)}\mathcal{M}(\Sigma)$. 
\end{lem}
\textbf{Proof of Lemma \ref{deriv} : }We use Proposition \ref{mult} to determine how the default $D(\tau , \theta)$ behaves with the product of elements in $K(\Sigma ,-e^{\frac{i \pi \hbar}{2}})$. We work with one component $\sigma_{\chi}$ of the $\psi$-symbol at a time.
For $\gamma \in K(\Sigma , -1)$, we will note $E_{\gamma}=\chi (\Delta_{\gamma} +D_{\gamma})$, so that we can write $\sigma_{\chi}^{\gamma}(\tau , \theta , \hbar)=\sigma^{\gamma}(\tau , \theta ,0) +\hbar E_{\gamma}(\tau , \theta) + o (\hbar)$.
\newline Then, applying $\chi \in \hat{A}_{\Gamma}$ to Proposition \ref{mult} we have:
\begin{displaymath} \sigma_{\chi}^{\gamma \cdot \delta}(\tau , \theta , \hbar ) =\sigma_{\chi}^{\gamma}(\tau , \theta , \hbar) \sigma_{\chi}^{\delta}(\tau , \theta , \hbar ) + \frac{\hbar}{i} \underset{e}{\sum} \partial_{\tau_{e}}\sigma_{\chi}^{\gamma}(\tau , \theta , \hbar ) \partial_{\theta_{e}} \sigma_{\chi}^{\delta}(\tau , \theta , \hbar ) + o(\hbar) 
\end{displaymath}
We have $\sigma^{\gamma}_{\chi}(\tau ,\theta ,0)=f_{\gamma}(R_{\chi}(\tau, \theta)$. Recall that by the formula of Goldman given in Section \ref{sec:TQFT}, $f_{\gamma \cdot \delta}= f_{\gamma}f_{\delta} +\hbar \frac{\pi}{i} \lbrace f_{\gamma} , f_{\delta} \rbrace +o(\hbar)$. So, isolating terms of order 1 in $\hbar$, we get:
\begin{eqnarray*}\frac{\pi}{i}\lbrace f_{\gamma} , f_{\delta} \rbrace(R_{\chi}(\tau ,\theta)) + E_{\gamma \cdot \delta}(\tau , \theta)
\\
=E_{\gamma}(\tau , \theta )f_{\delta}(R_{\chi}(\tau , \theta )) +E_{\delta}(\tau , \theta) f_{\gamma}(R_{\chi}(\tau , \theta))+\frac{1}{i}\underset{e}{\sum} \partial_{\tau_{e}}f_{\gamma}(R_{\chi}(\tau ,\theta )) \partial_{\theta_{e}}f_{\delta}(R_{\chi}(\tau ,\theta))
\end{eqnarray*}
but $\lbrace f_{\gamma} , f_{\delta} \rbrace =\frac{1}{2 \pi} \underset{e}{\sum}\partial_{\tau_{e}}f_{\gamma}\partial_{\theta_{e}}f_{\delta} - \partial_{\tau_{e}}f_{\delta}\partial_{\theta_{e}}f_{\gamma}$. We deduce that 
\begin{displaymath}E_{\gamma \cdot \delta}
=E_{\gamma}\sigma_{\chi}^{\delta} +E_{\delta} \sigma_{\chi}^{\gamma}+\frac{1}{2 i}\underset{e}{\sum} \partial_{\tau_{e}}\sigma_{\chi}^{\gamma} \partial_{\theta_{e}}\sigma_{\chi}^{\delta} +\partial_{\theta_{e}}\sigma_{\chi}^{\gamma} \partial_{\tau_{e}}\sigma_{\chi}^{\delta}
\end{displaymath}
However, as for $\gamma,\delta \in K(\Sigma ,-1)$ we have $f_{\gamma \cdot \delta}=f_{\gamma}f_{\delta}$, and $\chi (\Delta_{\gamma})=\frac{1}{2 i}\underset{e}{\sum}\frac{\partial^{2} f_{\gamma}}{\partial \tau_{e} \partial \theta_{e}} \circ R_{\chi}$, the Leibniz rules implies that $\chi (\Delta_{\gamma})$ satisfies the same law of composition:
\begin{displaymath}\chi(\Delta_{\gamma \cdot \delta})
=\chi(\Delta_{\gamma})f_{\delta} +\chi(\Delta_{\delta}) f_{\gamma}+\frac{1}{2 i}\underset{e}{\sum} \partial_{\tau_{e}}f_{\gamma} \partial_{\theta_{e}}f_{\delta} +\partial_{\theta_{e}}f_{\gamma} \partial_{\tau_{e}}f_{\delta}
\end{displaymath}
This concludes the proof of Lemma \ref{deriv}: $\chi \circ D$ is a derivation.
\newline
\newline
\textbf{Proof of Lemma \ref{tang} : }It is well-known that $\mathcal{M}'(\Sigma)$ is an affine algebraic variety  those smooth points is the open dense subset $\mathcal{M}_{\textrm{irr}}'(\Sigma)$ (see \cite{Si09} for example). The point $R_{\chi}(\tau , \theta)$ is thus a smooth point of $\mathcal{M}'(\Sigma)$ for any $(\tau ,\theta) \in U \times \mathbb{R}/2 \pi \mathbb{Z}$.
\newline
Then the proof comes from elementary considerations of algebraic geometry: when $V$ is an affine algebraic variety, and $x$ a point of $V$, we put a structure of $\textrm{Reg}(V)$-module on $\mathbb{C}$ by defining $f \cdot \lambda =f(x) \lambda$. Then $\textrm{Der}_{x} (V,\mathbb{C})$ identifies with $T_x V =m_x /(m_x)^{2}$ the algebraic tangent space to $V$ at $x$ (where $m_x =\lbrace f \ / \ f(x)=0 \rbrace$), and the algebraic tangent space at a smooth point is the same of as the tangent space of $V$ at $x$ in the sense of differential manifolds. As the affine variety $\mathcal{M}'(\Sigma)$ is smooth on the image of $R_{\chi}$, by this general property, derivations of $\textrm{Reg}(\mathcal{M}(\Sigma))$ can be viewed as vectors of the tangent space.  As $(\tau ,\theta) \mapsto R_{\chi}(\tau ,\theta)$ is a parametrization of $\mathcal{M}(\Sigma)$, the vector fields $(\partial \tau_e ,\partial \theta_e)$ give a basis of the tangent space $T_{R_{\chi}(\theta ,\tau)}\mathcal{M}(\Sigma)$ for each $(\tau ,\theta)$.
\newline
\newline
\textbf{Proof of Proposition \ref{default}} from the Lemmas \ref{deriv} et \ref{tang}. Combining these two lemmas allow us to assert that $\chi(D(\tau , \theta ))$, viewed as a map $\textrm{Reg}(\mathcal{M}'(\Sigma)) \rightarrow \mathbb{C}$, is of the form  $f \rightarrow \mathcal{L}_{X} f (R_{\chi}(\tau ,\theta))$ for some $X \in T_{R_{\chi}(\tau ,\theta )}\mathcal{M}'(\Sigma)$ and we may write
\newline $X=\underset{e}{\sum} a_{e} \frac{\partial}{\partial_{\tau_{e}}} + b_{e} \frac{\partial}{\partial_{\theta_{e}}}$ for some coefficients $a_{e},b_{e} :\mathcal{M}(\Sigma) \rightarrow \mathbb{C}$. As the default is smooth, so are the coefficients $a_e$ and $b_e$.  
\newline We want to prove that these coefficients all vanish. To this end, we remind that we proved in Section \ref{sec:ex} that the default vanishes for the curves $C_e$ and $D_e$. Besides, we have the formula of Section \ref{sec:comput}:
\newline We have $\sigma^{C_e}(\tau ,\theta ,\hbar)=-2 \cos (\pi \tau_e)[0]$, so that $\chi(D_{C_e})(\tau ,\theta)=2 a_e \pi \sin(\pi \tau_e)$. As the default vanishes on $C_e$, we must have $a_e=0$.
\newline To show the vanishing of the $b_e$, we use the formulas for $D_e$:
\newline In the first kind of curve $D_e$, described in Section \ref{sec:ex}, we have $f_{D_e}(R_{\chi}(\tau , \theta))=\sigma_{\chi}^{D_e}(\tau , \theta ,0) =2  W(\pi \tau ,0) \cos(\theta_e)$ where $W$ does not vanish for $\tau \in U$.
\newline We know that default $D_{D_e}$ vanishes, so we have
\newline $\chi(D_{D_e}(\tau , \theta ))=b_e \frac{\partial}{\partial \theta_e}f_{D_e}(R_{\chi}(\tau , \theta))=-2 b_e \pi \sin (\theta_e ) W(\pi  \tau , 0)=0$. This yields $b_e=0$. 
\newline 
\newline In the second case, we have $f_{D_e}(R_{\chi}(\tau ,\theta ))=\sigma_{\chi}^{D_e}(\tau ,\theta , 0)= -2 J(\pi \tau ,0) \cos (2 \theta_e) -I(\pi \tau ,0)$ for the  functions $I$ and $J$ defined in Section \ref{sec:ex}, that are non-vanishing for $\tau \in U$. 
\newline Again as
\newline $\chi (D_{D_e}(\tau ,\theta))=b_e \frac{\partial}{\partial \theta_e} f_{D_e}(R_{\chi}(\tau ,\theta))= 4 \pi b_e \sin (2 \theta_e) J(\pi \tau ,0)$ vanishes, we must have $b_e=0$. It follows that the default $\gamma \mapsto D_{\gamma}$ is the zero derivation on $K(\Sigma ,-1) \mapsto A_{\Gamma}$, which is the last ingredient we needed to complete the proof of the Proposition \ref{default}.

\end{document}